\documentclass[11pt,notitlepage,twoside,a4paper]{amsart}
\usepackage{amsmath}
\usepackage{amsfonts}
\usepackage{amssymb}
\usepackage{amscd}
\usepackage{amsthm}
\usepackage{amsbsy}
\usepackage[french, british]{babel}

\usepackage{amsmath,amssymb,enumerate}

\usepackage{epsfig,fancyhdr,color}%,showkeys,amsmidx?Aris

\usepackage{amssymb}
\usepackage{amsmath,amsthm}
\usepackage{latexsym}
\usepackage{amscd}
\usepackage{psfrag}
\usepackage{graphicx}
\usepackage[latin1]{inputenc}
\usepackage[all]{xy}
\usepackage[mathcal]{eucal}

\definecolor{NoteColor}{rgb}{1,0,0}

\renewcommand{\textsc}{\textcolor{red}}

\newtheorem*{theorem 1}{\rm\bf Proposition 1}
\newtheorem*{theorem 2}{\rm\bf Proposition 2}

\theoremstyle{definition}

\theoremstyle{remark}

\def\interieur#1{\mathord{\mathop{\kern 0pt #1}\limits^\circ}}

\begin{document}
%------------------------------------------------------------------------------------------------------------------
\title[A glimpse into Thurston's work]{A glimpse into Thurston's work}

\author[Ken'ichi Ohshika and Athanase Papadopoulos]{Ken'ichi Ohshika and Athanase Papadopoulos}

 \address{Ken'ichi Ohshika, Department of Mathematics,
Gakushuin University,
Mejiro, Toshima-ku, Tokyo, Japan}
  \email{ohshika@math.gakushuin.ac.jp}

 \address{Athanase Papadopoulos, Institut de Recherche Math\'ematique Avanc\'ee\\ CNRS et Universit\'e de Strasbourg\\\small 7 rue Ren\'e
  Descartes - 67084 Strasbourg Cedex, France}
  \email{papadop@math.unistra.fr}

\date{}
\maketitle

  \begin{abstract}

We present an overview of some significant results of Thurston and their impact on mathematics. The final version of this article will appear as Chapter 1 of the book \emph{In the tradition of Thurston: Geometry and topology}, ed. K. Ohshika and A. Papadopoulos, Springer Verlag, 2020.

  \bigskip

\noindent {\bf Keywords} foliation, contact structure, symplectic structure, conformal geometry, holomorphic motion,  Kleinian groups,   circle packing,  automatic groups, tiling,  mapping class groups, Teichm\"uller space, fashion design, linkage,  Anti-de Sitter geometry, higher Teichm\"uller theory, Grothendieck--Thurston theory,  asymmetric metric, Schwarzian derivative, computer science,
 Ehrenpreis conjecture, transitional geometry, 3-manifolds, geometric structures, (G,X)-structures, Dehn surgery, hyperbolic geometry, Thurston norm, Smith conjecture, Cannon--Thurston map, discrete conformal geometry, discrete Riemann mapping theorem.
   
\medskip

\noindent {\bf AMS codes} 57N10, 57M50, 20F34, 20F65, 22E40, 30F20, 32G15, 30F60, 30F45, 37D40, 57M25, 53A40, 57D30, 58D05, 57A35, 00A30, 01A60, 20F10, 68Q70, 57M05, 57M07, 57Q15, 57D15,  58A10, 58F10, 65Y25.

   \end{abstract}

\medskip
\medskip

\tableofcontents
\section{Introduction}

In this chapter, we present an overview of some significant results of Thurston and their impact on mathematics.

The chapter consists of two parts. In the first part,  we review some works of Thurston,  grouped in topics. The choice of the topics reflects our own taste and our degree of knowledge. The choice of the order of these topics was almost random. Indeed, it is not clear whether a given topic is more important than another one, and there are interconnections and mutual influences between most of these topics.
Furthermore, it was not possible to follow a chronological order because Thurston was thinking about all these subjects simultaneously.

In the second part of this chapter, we report briefly on the proofs of some conjectures which were either  formulated by Thurston or whose solution depended in a crucial way on his work. We also discuss 
 a few topics whose development was  directly or intellectually influenced by ideas of Thurston.
 
 We have included at some  places remarks and quotations which give an idea of Thurston's approach to science in general and to the aesthetics of mathematics.   

Our exposition will certainly be too short at some places, for some readers who know little about Thurston's work on the topic discussed, and it will be redundant  for readers  familiar with this topic (and even more for the experts). We apologize in advance to both categories of readers. We have added here and there some historical notes, whenever we felt this was useful.  These notes will probably be beneficial to both groups of readers.

We would like to thank Vlad Sergiescu for his comments on the section on foliations.

 \section{On Thurston's works}
\subsection{Foliations and groups of homeomorphisms}

 The first time the word\index{foliation} ``foliation" was used in a mathematical sense (in its French version, \emph{feuilletage}) took place by the end of the 1940s by Georges Reeb and Charles Ehresmann (who was Reeb's advisor).\footnote{In their first papers on the subject, Ehresmann and Reeb used the expression \emph{\'el\'ements de contact de dimension $p$ complètement int\'egrables} (``completely integrable contact elements of dimension $p$.")}   Reeb, in his dissertation \cite{Reeb}, gave the first example of a foliation on the 3-sphere.\footnote{The question of the existence of a foliation on the 3-sphere was asked by Heinz Hopf in 1935, who certainly did not use the word ``foliation". The first example was given in a joint paper by Ehresmann and Reeb, but Ehresmann always attributed this construction to Reeb.}

  When Thurston came into the subject, examples and constructions of foliations on special manifolds were available. To describe the situation in short, one can say that in a lapse of time of  five years, he obtained all the general existence results that were hoped for. In this section, we briefly review his work on the subject. 
  
 Thurston's first published paper on foliations\footnote{With the exception of his PhD thesis, which was defended the same year and  which remained unpublished. The thesis, whose title is \emph{foliations of 3-manifolds which fiber over a surface.}
was submitted to the Swiss journal \emph{Inventiones}. 
 The referee asked for modifications;  Thurston did not comply and withdrew the paper. Haefliger,  in the collective Thurston memorial article \cite{Gabai-Memorial}, writes: ``The referee suggested
that the author should give more explanations. As
a consequence, Thurston, who was busy proving
more theorems, decided not to publish it."} is a short note titled \emph{Noncobordant foliations of $S^3$} \cite{Thurston1972}, which appeared  in 1972. In this note, Thurston proved that any closed 3-manifold carries a family of foliations whose Godbillon--Vey invariant takes all possible real values.  This invariant (an element of the real 3-cohomology of the manifold) was discovered the year before, by C. Godbillon and J. Vey,  who came across it by manipulating differential forms. The question of  whether there exist foliations with non-zero Godbillon--Vey\index{Godbillon--Vey invariant} invariant was soon raised, together with the problem of giving it a geometrical interpretation. Thurston's result closed both problems. It is interesting to note that Thurston's existence proof has a strong hyperbolic geometry flavor: he constructed a family of foliations of the 3-sphere which depend on convex polygons in the hyperbolic plane whose area is equal to the Godbillon--Vey invariant of the foliation.  Besides hyperbolic geometry, we can find in this proof another ingredient which was soon to become fundamental in Thurston's work, namely, the notion of singular hyperbolic surface. In his paper on foliations, Thurston described these surfaces as ``surfaces having a number of isolated corners, with metrics of constant negative curvature everywhere else."
 
  Thurston formulated his result on the Godbillon--Vey invariant as the surjectivity of a certain homomorphism from the group of cobordism classes of foliations onto $\mathbb{R}$. As a corollary, he proved the existence of an uncountable family of non-cobordant foliations on $S^3$.
The precise  results are stated in terms of Haefliger's classifying\index{Haefliger classifying space}\index{Haefliger structure} spaces of $\Gamma$-structures. These  objects, also called Haefliger structures, were introduced by Andr\'e Haefliger in his thesis published in 1958 \cite{Haefliger1958}, as a generalization of the notion of foliation.\footnote{\emph{A posteriori}, it is interesting to read Richard Palais' comments on this notion, in his Mathscinet review of  Haefliger's paper:  ``The first four chapters of the paper are concerned with an extreme, Bourbaki-like generalization of the notion of foliation. After some twenty-five pages and several hundred preliminary definitions, the reader finds that a foliation of $X$ is to be an element of the zeroth cohomology space of $X$ with coefficients in a certain sheaf of groupoids.  While such generalization has its place and may in fact prove useful in the future, it seems unfortunate to the reviewer that the author has so materially reduced the accessibility of the results, mentioned above, of Chapter V, by couching them in a ponderous formalism that will undoubtedly discourage many otherwise interested readers." In fact, the notion that Haefliger introduced turned out to be of paramount importance.} In short, Haefliger structures, which may be interpreted as singular $C^r$ foliations, are  $\mathbb R^k$-bundles over $n$-dimensional manifolds equipped with foliations transverse
to the fibers. (Such a Haefliger structure is said to be of codimension $k$). A natural example of a Haefliger structure is the normal bundle to a foliation. Haefliger structures are the natural setting for the theories of classifying spaces and of characteristic classes of foliations, and  Haefliger's theory reduces the question of the existence of certain classes of foliations to that of certain maps between manifolds and classifying spaces. A consequence of
 Thurston's work on Haefliger structures is that  in some sense (up to a natural condition on normal bundles) the class of $\Gamma$-structures is not different from that of foliations. 
   
 Two years later, Thurston published three papers in which he proved a series of breakthrough results on foliations. In the paper \cite{thurston-cmh}, titled \emph{The theory of foliations of codimension greater than one}, working in the setting of codimension-$k$ 
Haefliger structures for $k>1$, he showed the existence of a large class of completely integrable  plane fields on manifolds which led to the construction of new classes of foliations. In particular, he obtained that any plane field of codimension greater than one is homotopic to a completely integrable $C^0$ field. He also proved that for any $n\geq 3$ and any $1< k\leq n/2$, if the sphere $S^n$ carries a $k$-plane field then it carries also a $C^\infty$ foliation of dimension $k$.  These results are wide generalizations of existence results obtained in special cases by Reeb, Tamura, Lawson, Phillips, Haefliger, and others.

The second paper by Thurston published in 1974 \cite{Thurston-stability} is  titled \emph{A generalization of the Reeb stability theorem}. The theorem, obtained by Reeb in his thesis \cite{Reeb} which we already mentioned, says that if a transversely oriented codimension-1 foliation on a compact manifold  has a two-sided compact leaf with finite fundamental group, then 
 all the leaves of this foliation are diffeomorphic.  Thurston obtained a much stronger result under the hypothesis that the foliation is of class $C^1$, namely, he proved that one can replace the hypothesis that the compact leaf has finite fundamental group by the one saying that the first real cohomology group  of the leaf is zero. At the same time, he showed that, under the same conditions, the leaves of this foliation are the fibers of a fibration of the manifold over the circle or the interval. Furthermore, he showed by an example that this result does not hold in the $C^0$ case. Thurston approached the problem by studying the linear holonomy around the compact leaf, and he gave an interpretation  of the result in terms of the linearity properties at a fixed point for a topological group acting continuously in the $C^1$ topology, as a group of $C^1$ diffeomorphisms of a manifold.

The relation between foliations and groups of diffeomorphisms is announced  in the title of the third paper published in the same year \cite{Foliations-and-groups}:  \emph{Foliations and groups of diffeomorphisms}. In this paper,  Thurston studied higher codimension Haefliger structures in relation with groups of diffeomorphisms of arbitrary manifolds, generalizing a relation discovered by John Mather between 
 the group of compacty supported diffeomorphisms of the real line and 
framed codimension-one Haefliger structures. Using the techniques of classifying spaces, Thurston proved that any two $C^\infty$ foliations of a manifold arising from nonsingular vector fields are homotopic as Haefliger structures if and only if their normal bundles are isotopic. At the same time, he announced a result giving a precise relationship between the classifying space of codimension-$k$ Haefliger structures and that of the diffeomorphism group of compact manifolds of dimension $k$ as a discrete group. He used this result to prove that the  diffeomorphism group of a compact manifold is perfect (that is, equal to its commutator subgroup), generalizing a result obtained by Mather  in the case of 1-dimensional manifolds.  

One may recall here that at the time Thurston was working on these topics, the study of the algebraic structure of groups of diffeomorphisms\index{group of diffeomorphisms of a manifold} and homeomorphisms of compact manifolds was a very active subject of research, involving mathematicians such as John Mather, David Epstein,  Michel Herman, Jean Cerf and others.
  
In the same year, Thurston gave a talk at the Vancouver ICM (1974) whose title was \emph{On the construction and classification of foliations} \cite{Thurston-Vancouver}. In the paper  published in the proceedings of the congress, Thurston reviewed some of the major results he had obtained and he announced the results of his forthcoming paper \cite{T-existence}, \emph{Existence of codimension-one foliations}, which appeared in 1976. In this paper, he proved the existence of a $C^\infty$ codimension-one foliation on any closed manifold whose Euler characteristic is zero. This result may be contrasted with a result of Haefliger \cite{Haefliger1958} stating that there is no codimension-one real-analytic foliation on a sphere of any dimension. This displays a striking difference between the  $C^{\infty}$ and the real-analytic cases.
In the same paper, Thurston proved that on any closed manifold without boundary,  every hyperplane field is homotopic to the tangent plane field of a $C^\infty$-foliation.  

Naturally, the ICM paper, which is only 3 pages long, is written in the pure Thurston style, warm, unconventional and appealing to the reader's imagination. It starts with the following: 
\begin{quote}\small
Given a large supply of some sort of fabric, what kinds of manifolds can be made
from it, in a way that the patterns match up along the seams? This is a very general
question, which has been studied by diverse means in differential topology and differential
geometry.
\end{quote}
 
 It is also not surprising that the definition of a foliation that Thurston gives in this paper is informal and unusual: 
 \begin{quote}
 A foliation is a manifold made out of a striped fabric---with infinitely thin stripes, having no space between them. The complete stripes, or ``leaves,"  of the foliation are submanifolds; if the leaves have codimension $k$, the foliation is called a 
codimension-$k$ foliation.
\end{quote}

ICM talks are intended for a general audience, but very few mathematicians were able, like Thurston was, to describe the objects they were studying in simple words, avoiding notation and formulae.

 Thurston's paper \cite{T-Norm} which appeared in 1986 and which is titled \emph{A norm for the homology of 3-manifolds},\index{Thurston norm} is the foundational paper on the so-called Thurston norm,\index{Thurston norm} and it also contains results and conjectures on foliations of 3-manifolds. The results include a classification of codimension-1 foliations without holonomy that are transverse to the boundary, in terms of the top-dimensional faces of the unit ball of Thurston's norm on homology. We shall talk about hhis in \S \ref{s:norm} below. Thurston showed that if such a foliation has no Reeb component, then any compact leaf is norm-minimizing in its homology class. A converse was obtained by D. Gabai soon after,  who also obtained a general existence result for codimension-1 transversely oriented foliations transverse to the boundary with no Reeb components \cite{Gabai1982}. In the same paper, Gabai proved several conjectures of Thurston.
 
 In 1976, Thurston published a paper with Joseph Plante on the growth\index{growth! of germs of diffeomorphisms} of germs of diffeomorphisms\index{germs of diffeomorphisms} \cite{PT}. The study of such germs  was motivated by the theory of foliations (the germs appear in the holonomy groups of foliations).
The question that Plante and Thurston studied was motivated by the work of Milnor  \cite{Milnor1968} who introduced the notion of growth of a finitely generated group.\index{growth! of finitely generated groups} It was conjectured that a group has polynomial growth if and only if it contains a nilpotent subgroup of finite index. The conjecture was proved by Gromov in a paper which was a major breakthrough \cite{Gromov1981}.
In their article, Plante and Thurston showed that the conjecture is true for germs of diffeomorphisms, and they gave applications of this result to foliations. One of the consequences they obtained is that if a compact manifold with fundamental group of polynomial growth carries a transversely oriented and real-analytic codimension-1 foliation, then its first real homology group is nontrivial.
This generalized a result of Haefliger \cite{Haefliger1958}.
 
 One should also talk about Thurston's mostly unpublished work on extending the theory of foliations to that of laminations and essential laminations in 3-manifolds which was developed by several authors, following his ideas.  
 
  We mention Thurston's manuscript \cite{Taut} on taut foliations,\index{taut foliation}\index{foliation!taut} that is, codimension-1 foliations for which there exists an embedded closed curve that is transverse to the foliation and intersects every leaf.\footnote{The terminology \emph{taut foliation} is also used for a foliation such that the ambient manifold admits a Riemannian metric for which the leaves are minimal surfaces. In his paper \cite{Sullivan1979}. Sullivan showed that for a $C^2$ taut foliation (in the above sense) on a closed orientable 3-manifold, there exists a Riemannian metric on the ambient manifold that is taut in this Riemannian sense.}
 In this paper, Thurston associates to every transversely orientable taut foliation on a closed atoroidal 3-manifold a faithful homomorphism from the fundamental group of the manifold onto the group of orientation-preserving homeomorphisms of the circle which preserves a pair of dense invariant laminations of the circle (in an appropriate sense), and which is universal in some precise sense. This theory gave rise to developments by several authors, see in particular the paper  \cite{Calegari-Dunfield} by Calegari and Dunfield in which the authors give a new proof of Thurston's result and  where Thurston's ideas are made more precise. In the same paper, the authors show that there are other classes of essential foliations and laminations than those which were considered by Thurston that give rise to faithful actions on the circle.  
 We also refer the reader to the exposition in Chapter 10 by Baik and Kim in the present volume \cite{Baik-Kim}. See also  the book \cite{Calegari} by Calegari on foliations which contains several sections explaining Thurston's homomorphism to the group of orientation-preserving homeomorphisms of the circle. As a matter of fact, this books is a valuable reference for many aspects of Thurston's theory of foliations and laminations,  including his work on the cohomology of the group of orientation-preserving homeomorphisms of the circle, his stability result for the group of orientation-preserving homeomorphisms of the interval, his construction of foliations on 3-manifolds using a triangulation of the manifold, the Thurston norm in relation with foliations, and several other topics.

Talking about Thurston's work on foliations, one should also mention measured foliations on surfaces and his construction of the space of measured foliations, a space equipped with a topology which makes it homeomorphic to a Euclidean space of the same dimension as Teichm\"uller space. The space of measured foliations became a central object in low-dimensional geometry and topology. We review this in \S \ref{s:surfaces} below.
  
 \subsection{Contact and symplectic geometry}  \label{s:symplectic}

A contact structure on a differentiable manifold is a field of hyperplanes in the tangent bundle satisfying a ``complete non-integrability" condition that makes it, locally,  unrealizable as a hyperplane field tangent to a foliation. (Note that the non-integrability property of these hyperplanes  is in contrast with the dimension-1 case: vector fields are locally always integrable.) It is easy, although not trivial, to produce examples of contact structures. Several examples are discussed in Thurston's book on the geometry and topology of 3-manifolds \cite{Thurston1}. 

Contact geometry,\index{contact structure}\index{contact geometry}\index{symplectic geometry} like symplectic geometry which we shall discuss below, originates in classical mechanics, and it has applications in geometric optics, thermodynamics and other domains of physics. In fact, this notion can be traced back to the work of Gaston Darboux. One of his results is often quoted, viz. the fact that a contact structure is always locally equivalent to a standard contact structure \cite{Darboux-sur}. 

The usual definition of a contact structure is algebraic, formulated in terms of differential forms. In fact, a foliation and a contact structure are both defined locally by a differential form $\alpha$, but in the case of a foliation, $\alpha$ satisfies $\alpha\wedge d\alpha\equiv 0$ whereas in the case of a contact function, it satisfies $\alpha\wedge d\alpha \not=0$  at every point (and if the relation is replaced by $\alpha\wedge d\alpha >0$ with respect to a given orientation, we say that we have a ``positive" contact form). The problem of classifying contact structures on manifolds arose naturally. 
Thurston writes in his monograph \cite[p. 168]{Thurston1}: ``[contact structures] give an interesting example of a widely occurring pattern for manifolds that is hard to see until your mind and eyes have been attuned."  Several pages of his book \cite{Thurston1} are dedicated to the effort of sharing with the reader an intuitive picture of contact structures.  On p. 172 of this book, to give a physical sense of the contact structure on the tangent circle bundle of a surface, he uses models related to ice skating and bicycling,  dedicating several paragraphs to these images.

Thurston obtained a number of important results on contact structures.     His first paper on this subject is a joint paper with H. E. Winkelnkemper \cite{TW}, titled \emph{On the existence of contact forms},\index{contact form} published in 1975. It contains a very short proof of a result, which was already obtained by Robert Lutz and by Jean Martinet in 1971 \cite{Lutz, Martinet},
 saying that every closed orientable 3-manifold carries a contact structure. Thurston and Winkelnkemper  deduced this result from a classical result, namely, the so-called ``open-book decomposition theorem" of Alexander \cite{Alexander} (1923).  
  
 Even though, in some sense, a contact structure is the complete opposite of a foliation,   the two subjects are related. With Yakov Eliashberg, Thurston introduced the notion of  \emph{confoliation}\index{confoliation} in dimension 3, and he developed a theory which gives a hybrid setting for codimension-1 foliations and contact structures on 3-manifolds \cite{ET1, ET2}. A confoliation in this sense interpolates between a codimension-1 foliation and a contact structure. The techniques that Eliashberg and Thurston developed allowed them to prove that any $C^2$ codimension-1 foliation on a 3-manifold, except for the product foliation $S^1\times S^2$, can be approximated in the $C^0$ sense by positive contact structures. Confoliations appear in a crucial manner in the proof of this result, since the main step  consists in the modification of the plane field tangent to a foliation so that it becomes integrable (tangent to a foliation) in some part of the manifold and a positive contact structure in the complement. It is interesting to note that at the same time he was developing confoliations, Thurston developed a theory of foliations  of three-manifolds that are hybrids of fibrations over the circle and foliated circle bundles over surfaces, see his 1997 preprint
 \cite{Taut}.

  Contact structures\index{contact structure} are defined on odd-dimensional manifolds, and their analogues on even-dimensional manifolds are symplectic structures.\index{symplectic structure}
  
  This brings us to Thurston's work on symplectic geometry.\index{symplectic geometry}

When Thurston started working in this field,  the questions of the existence of symplectic structures on closed manifolds and  that of K\"ahler metrics on symplectic manifolds were the main problems. 
 In 1976,  Thurston gave the first examples of compact symplectic manifolds that do not admit any K\"ahler metric.\index{K\"ahler metric} He presented his examples in a short note titled \emph{Some simple examples of symplectic manifolds} \cite{T-AMS1976}. The examples were later called Kodaira--Thurston manifolds, since  it turned out that the  manifolds described by Thurston were already known to Kodaira (who used them for other purposes). At the same time, Thurston gave a counter-example to a claim made by Heinrich Guggenheimer, saying that the odd-dimensional Betti numbers of symplectic manifolds are necessarily even. The examples that Thurston gave have odd first Betti numbers. The odd-dimensional Betti numbers of K\"ahler manifolds are all even. After Thurston's discovery, the question of characterizing the symplectic manifolds which admit no K\"ahler structure became a very active research field (works of Robert Gompf, Dusa McDuff, etc.).
  
One may also mention here a result of Thurston on  volume-preserving diffeomorphisms\index{volume-preserving diffeomorphism} of differentiable manifolds. This theory  is related to symplectic geometry, since a symplectomorphism   of a $2n$-dimensional manifold preserves the volume form obtained as the $n$-th  power of the symplectic form.   In a preprint titled \emph{On the structure of the group of volume preserving diffeomorphisms}, first circulated in 1972  \cite{Volume72}, Thurston proved that the group of volume-preserving diffeomorphisms of a manifold is perfect provided the first homology group of the manifold is zero, and he introduced at the same time a certain number of ideas that became later very useful in symplectic geometry. 
 Although Thurston's preprint remained unpublished, 
the techniques it contains and the questions it raises had a profound impact in symplectic geometry (works of Augustin Banyaga, of Dusa McDuff, etc.).
Banyaga,  in his paper \cite{Banyaga1} and in his book \cite[p. 125ff]{Banyaga}, developed many ideas of Thurston and gave detailed proofs of  several of his results in symplectic geometry. 
              \subsection{1-dimensional dynamics}

      The first published work by Thurston on dynamics is his paper with Milnor \emph{On iterated maps of the interval}, which appeared in 1988  \cite{MT}. An early version of the paper,  containing more material,  was circulated in 1977.\footnote{Leo Jonker, in his Mathscinet review of this paper writes: ``If there were a prize for the paper most widely circulated and cited before its publication, this would surely be a strong contender. An early handwritten version of parts of it was in the reviewer's possession as long ago as 1977."} Some results stated as a conjecture in the preprint version became theorems in the published version.

   In this paper,  Milnor and Thurston
studied the dynamics of continuous piecewise (strictly) monotone maps of the interval, to which they associated a certain number of naturally and very simply defined invariants.   These invariants became at the basis of kneading theory, an important element in the theory of   dynamics\index{dynamics of the interval} of unimodal maps\index{unimodal map} of the interval.  Let us recall some of the notions they introduced. 

Given a map $f$ of an interval $I$, a \emph{lap} of $f$ is a maximal sub-interval of $I$ on which $f$ is monotone. This leads to the notion of  \emph{lap number} $\ell=\ell(f)$ of $f$. Milnor and Thurston studied  the \emph{growth}\index{growth!of lap number} 
of the lap number of the iterates of $f$, that is, the limit $\lim_{k\to\infty} \ell(f^k)^{1/k}$. By a theorem of  Misiurewicz and Szlenk, this limit is equal to the topological entropy of $f$. Milnor and Thurston introduced an invariantly defined ``formal coordinate function" $\theta(x)$ which is given for $x$ in $I$, as a formal power series   $\sum \theta_k(x)t^k$, where if $f^k(x)$ belongs to the interior of the $j$-th lap $I_j$, the coefficient $\theta_k(x)$ is the formal symbol $I_j$ multiplied by $\pm 1$ or $0$ according as to whether $f^k$ is increasing or decreasing, or has a turning point at $x$.
  This led them to a basic invariant called the \emph{kneading matrix}\index{kneading matrix} of $f$, an $(\ell +1)\times \ell$ matrix with entries in the ring $\mathbb{Z}[[t]]$ of integer formal power series, with its associated \emph{kneading determinant}, a\index{kneading determinant} power series with odd integer coefficients, $D(t)= 1+D_1(t)+D_2(t)+\ldots$. There is a close relation between the kneading determinant and the behavior of the periodic points of the map. In the simplest case where $f$ has only one turning point (which is the critical point of the map), the coefficients of $D(t)$ are either $+1$, $-1$ or $0$ according to whether the iterate $f^{k+1}$ has a local minimum or a local maximum at $c$. In the same paper,
Milnor and Thurston gave a method for computing the sequence of lap numbers of the iterates of $f$ in terms of the kneading matrix.  
 They  studied the convergence of the kneading determinant, showing for example that for  $s>1$, this power series is holomorphic in the unit disc, and has a smallest zero at $t=1/s$  where $s=s(f)$  is the growth number of the map. Under the same hypothesis ($s>1$), they showed that $f$ is topologically semi-conjugate to a piecewise linear map having slope $\pm s$ everywhere. The Artin--Mazue zeta function encodes the periodic orbits of $f$. 
 
 Milnor and Thurston used methods of Julia and Fatou, before these methods found their place  in the revival of holomorphic dynamics that took place a few years later. They proved what they call their \emph{main theorem}, which allows a computation of the Artin--Mazur zeta function in terms of the kneading determinant. They gave several applications of their theory. 
 
 A particularly important class of examples of maps to which the Milnor--Thurston applies is the one of maps of lap number two, \emph{unimodal real maps}.\index{unimodal map}  A typical family of such maps is the family of quadratic polynomials $x\mapsto x^2+c$. Each map in this family has a unique critical point, and the kneading sequence\index{kneading sequence} describes the location of the sequence of images of this critical point, to the left or right of this critical point. 
 For the family of quadratic polynomials, Milnor and Thurston gave a characterization of power series that can occur as a kneading determinant, they discussed continuity properties of the growth number $s(f)$, and they obtained a monotonicity result for the entropy.
Furthermore, the paper contains several algorithms to compute the entropy of a piecewise monotone map.

  This is now the occasion for us to quote Milnor from the preface and the dedication to Thurston that he wrote, in Volume VI of his \emph{Collected Works} \cite{Milnor-collected-VI}.\footnote{The volume is dedicated to Thurston.} In the preface, Milnor writes: ``I was introduced to Dynamical Systems by Bill Thurston in the late 1970s and found the field so engaging that it was hard to escape from."
  In the dedication, Milnor writes: 
\begin{quote}\small
My interactions with Bill followed a consistent pattern. He would propose a mathematical statement which I found amazing, but extremely unlikely. However, the evidence would accumulate, and sooner or later I would have to concede that he was completely right. My introduction to the field of dynamics proceeded in exactly this way.
  Bill had been intrigued by the work of Robert May in theoretical ecology. May had proposed that the population of some insect species in successive years behaved in a chaotic way, which could be described by a very simple mathematical model, in which next years population is expressed as a universal modal function of this years population. Bill developed this idea by constructing symbol sequences associated with unimodal maps. He claimed that many quite different looking one-parameter families of unimodal maps would give rise to the same family of symbol sequences. I didn't believe a word of this, but couldn't find a counter-example. Eventually, I was convinced, and we collaborated on the paper ``On iterated maps of the interval."

  \end{quote}

    The notions that Milnor and Thurston introduced in their paper remain until now part of the most important tools for the study of the dynamics of maps of the interval. Their paper  continues to be a  source of inspiration for the works done in this field. A large volume of literature is devoted to the generalization of their results for unimodal maps to maps with a larger number of laps. Furthermore, kneading sequence theory, as a way of encoding combinatorial information, was applied in the study of complex dynamics, by Milnor and others.
    
           Thurston's last published paper \cite{Thurston2014}\footnote{The paper was published posthumously in 2014.} is on dynamics.  He wrote it before his death, a period where, according to Milnor, ``Bill was entering a period of renewed creativity, full of ideas and eager to communicate them." \cite[p. {\sc ix}]{Milnor-collected-VI} The paper is titled \emph{Entropy in dimension one}, and it contains new ideas and results in this field. One of the results that Thurston obtained is a characterization of positive numbers that arise as the topological entropy of postcritically finite self-maps of the interval. Precisely he proved that these are exactly the numbers $h$ such that $\exp(h)$ is an algebraic integer that is at least equal to the absolute value of any conjugate of $\exp(h)$. He also showed that the map can be chosen to be a polynomial whose critical points are all in the open interval $(0,1)$.
                            At the same time, the paper makes it clear what are the phenomena of 1-dimensional dynamics that are relevant for entropy.
                            
                          Thurston used in this paper a number of ideas and notions from his previous works:  the central role played by postcritically finite maps,  train tracks for graphs together with train track maps and the operations of zipping and splitting of train tracks (ideas originating in his theory of surface dynamics, under the version adapted by Bestvina and Handel in their study of outer automorphism groups of free groups),  a generalization of the notion of pseudo-Anosov maps, Perron--Frobenius matrices and Pisot and Salem numbers, notions that appear in Thurston's theory of surface automorphisms. 
    
    \subsection{The topology of 3-manifolds}
              \label{3-man}
             We start with a few words on the pre-Thurston era. 
              
In the 1960s, a new direction of research in 3-manifolds was started by Haken  and Waldhausen \cite{Ha, Wald}.
 Among the objects of their research is the class of compact irreducible 3-manifolds containing incompressible surfaces. These manifolds were called by Waldhausen sufficiently large 3-manifolds;\index{sufficiently large 3-manifold} now they are called Haken manifolds.\index{Haken manifold}
In particular, Waldhausen proved that any homotopy equivalence between two closed Haken manifolds is homotopic to a homeomorphism.
In the 1970s, Jaco, Shalen and Johannson developed a theory of decomposing Haken 3-manifolds along incompressible tori and annuli, setting the basis of a theory now called Jaco--Shalen--Johannson theory \cite{JS, Joh}.\index{Jaco--Shalen--Johannson theory}
They showed that any Haken manifold can be decomposed along (a possibly empty) union of  disjoint incompressible tori  in such a way that each piece is either a Seifert fibered manifold or an atoroidal manifold, i.e. a 3-manifold which contains no non-peripheral immersed incompressible tori.
(Here an immersed surface is said to be incompressible when the map induces a monomorphism between the fundamental groups.)

Through the work of Andreev, Riley and J{\o}rgensen, Thurston already noticed that there are many 3-manifolds that admit complete hyperbolic metrics.
He considered that this should be the case in much more generality.
Indeed  he proved a ``uniformization theorem for Haken manifolds,"  stating that any atoroidal Haken manifold which is closed or which has torus boundaries carries a complete hyperbolic metric of finite volume.\index{uniformization theorem!Haken manifolds}
His proof of this theorem is very intricate and long.
The argument is divided, at a large scale, into two cases: the first is when the manifold is not a surface bundle over the circle and the second is when it is.
The first case is proved by induction involving Maskit's combination theorem.
The second case is proved using his own theorem called the ``double limit theorem,"\index{double limit theorem} which is itself an important contribution to the theory of Kleinian groups, and which we describe in Section \ref{Kleinian}.
One of the remarkable consequences of the uniformization theorem is the resolution of the Smith conjecture, which we shall review in Section \ref{Smith}.

Thurston conjectured that the same  kind of uniformization theorem should hold for all closed 3-manifolds, and not only for Haken ones.
This conjecture was formulated in the form of a ``geometrization conjecture"\index{geometrization conjecture} which includes the Poincar\'{e} conjecture\index{Poincar\'e conjecture} as a very special case.\index{geometrization theorem}
The geometric structures\index{geometric structure} to which Thurston referred are  locally homogeneous metrics.\index{geometric structure}
He gave the list of eight kinds of three-dimensional geometric structures. Six among them can be carried only by Seifert fibered manifolds.
The two remaining ones are the hyperbolic geometry and the solvable geometry. Only torus bundles over the circle can carry a solvable geometry.
The geometrization conjecture\index{geometrization theorem} says that every compact 3-manifold  is decomposed along incompressible tori  into  3-manifolds having geometric structures.
In the case of a homotopy sphere, this is equivalent to the Poincar\'{e} conjecture.\index{Poincar\'e conjecture}

At the end of his the expository paper \cite{ThB}, Thurston gave a list of problems on 3-manifolds and Kleinian groups.
The problems on 3-manifolds contain the above-mentioned geometrization conjecture (and in particular  the Poincar\'{e} conjecture),\index{Poincar\'e conjecture} and the virtual-Haken conjecture,\index{virtual Haken conjecture} which says that every closed irreducible 3-manifold has a finite cover which is Haken. This question  was first posed by Waldhausen  (see \cite{Waldpr}).\index{virtual Haken conjecture}
The list also contains  quite a new and unexpected conjecture, now called the ``virtual fibering conjecture,"\index{virtual fibering conjecture}  which was proved by Agol more than a quater of a century later.

The list  as a whole has been the driving force of  all research in 3-manifolds and Kleinian groups for more than 30 years after its appearance.

  \subsection{$(G,X)$-structures and Geometric structures } 
 In the Erlangen programme, Klein proposed a new way of  thinking geometry.
 According to him, geometry consists of a base space and a group acting on it.  
 Although he did not think of general manifolds (the notion of manifold did not exist yet), we can regard his work as the origin of $(G,X)$-structures\index{G,X@$(G,X)$-structure}.
 A formal definition of a $(G,X)$-structure first appeared in Ehresmann's work. For a geometrico-historical exposition, we refer the reader to Goldman's article \cite{Goldman}.

Given a space $X$ and a group $G$ acting on $X$ by homeomorphisms, a manifold is said to have a $(G,X)$-structure when it is equipped with an atlas whose charts are maps into $X$ with transition maps being restrictions of elements of $G$.
 Thurston made this notion central in low-dimensional topology by giving many important and interesting examples of $(G,X)$-structures on manifolds. By his work,  this notion moved to the forefront of research. 
 
Geometric structures\index{geometric structure} are typical $(G,X)$ structures, where $X$ is a homogeneous space and $G$ is its isometry group.
As we mention in Section \ref{3-man}, Thurston showed that there are eight kinds of geometric structures in dimension $3$, and conjectured that every compact 3-manifold can be decomposed along incompressible tori into submanifolds having geometric structures.

Geometric structures are Riemannian structures, i.e., the stabilizer of $G$ at a point in $X$ is compact.
Thurston also studied non-Riemannian $(G,X)$-structures, above all  complex projective structures on surfaces.\index{complex projective structure}
This is the case where $X=\mathbb CP^2$ and $G=\mathrm{PSL}(2,\mathbb C)$.
The space of (marked) complex projective structures\index{complex projective structure} modulo isotopies on a closed surface $\Sigma$ had been studied from the viewpoint of complex analysis using Schwarzian derivatives by Bers, Kra and Marden among others.
Thurston gave a new parametrization for this space based on a more geometric approach, which has the form of a homeomorphism between this space and the product $\mathcal{T}(\Sigma) \times \mathcal{ML}(\Sigma)$, where $\mathcal T(\Sigma)$ denotes the Teichm\"{u}ller space of $\Sigma$ and $\mathcal{ML}(\Sigma)$ the space of measured laminations on $\Sigma$.
This more geometric approach to the space of complex projective structures opened up a new flourishing field, which should be called a topological study of projective structures. The reader is referred to \S \ref{s:grafting} of the present chapter.

 \subsection{Geometrization of cone-manifolds}\label{geometrization}
 After proving the uniformization theorem for Haken manifolds, Thurston tackled the general geometrization problem by a quite different approach.
 Since non-Haken manifolds do not contain  incompressible surfaces, there is no way to cut them  into simpler ones.
 Instead, Thurston introduced the technique of deforming of the structure of a cone manifold by increasing its cone angle.
 For instance, to prove that a non-Haken atoroidal 3-manifold $M$ has a hyperbolic structure, we would take a hyperbolic knot $K$ in $M$ (i.e. a knot $K$ such that $M \setminus K$ has a complete hyperbolic metric, which is guaranteed to exist by virtue of the uniformization theorem for Haken manifolds combined with Myers' theorem \cite{Myers}), and consider a deformation of the complete hyperbolic structure on $M\setminus K$ to a cone hyperbolic structure whose singular locus is the knot $K$, with cone angle $\theta$.
 If we could deform the cone-hyperbolic structure until the cone angle becomes $2\pi$, then we would be able to show that $M$ is hyperbolic.
 Of course, this strategy should break down in general, for $K$ may not be isotopic to a closed geodesic in a hyperbolic 3-manifold.
 
 What Thurston really proved can be expressed as follows, if we only consider the case where the singularity is a knot and there is no incompressible torus disjoint from the singularity. 
 Suppose that $M$  is a closed irreducible 3-manifold containing a prime knot $K$, and that we are given an angle $\theta=\pi/n$ on $K$.  
 We consider a 3-orbifold $(M,K(\theta))$ whose underlying space is $M$ and whose singular set is $K$ with cone angle $\theta$.
 Thurston proved that in this situation $(M,K(\theta))$ has some (possibly empty) disjoint incompressible Euclidean 2-sub-orbifolds which decompose $M$ into geometric 3-orbifolds.
 For simplicity, we consider the case where $M \setminus K$ is atoroidal.
 To prove the theorem (in this case), Thurston considered a deformation of the hyperbolic cone structure by increasing the cone angle on $K$, starting from the complete hyperbolic metric on $M$, which is regarded as the cone angle $0$.
% that for any closed $3$-manifold $M$ and a hyperbolic knot $K$ in $M$, we can deform the hyperbolic cone structure increasing the cone angle until the structure degenerates.
If the angle reaches $\theta$ without degeneration, then $(M, K(\theta))$ is a hyperbolic orbifold.
Thurston showed that if the degeneration occurs, then $(M, K(\theta))$ admits either a decomposition along an incompressible Euclidean 2-sub-manifold or  a geometric structure other than the hyperbolic one.
To prove the last step, in a special case, Thurston made use of  the Ricci flow and of Hamilton's theorem \cite{Ham} to get a spherical structure in the limit.
Thurston's geometrization conjecture,\index{geometrization conjecture} including the Poincar\'{e} conjecture,\index{Poincar\'e conjecture}  was solved later by Perelman using precisely these Ricci flows,\index{Ricci flow} based on the idea of Hamilton.
It is noteworthy that Thurston already noticed the usefulness of Ricci flows back in the 1980s.

The geometrization theorem of 3-orbifolds\index{geometrization theorem!3-orbifolds} implies that if a closed prime 3-manifold has a finite group action with one-dimensional fixed point set then it has  a geometric structure which is invariant under the action.
This is (quite a huge) generalization of the Smith conjecture.
 
\subsection{Dehn surgery}
Besides the Smith conjecture which we mentioned above, Thurston's work had a great impact on  knot theory,\index{knot theory} through his theory of hyperbolic Dehn surgery.\index{hyperbolic Dehn surgery}
Dehn surgery is a classical tool in knot theory.
\index{Dehn surgery}
The definition is simple: performing a Dehn surgery along a knot means that we take a tubular neighborhood of the knot and glue back the removed solid torus in such a way that the boundary of the meridian is in a homotopy class (called the meridian slope) on the boundary of the tubular neighborhood different from the original.
In this way, we get a new 3-manifold.
Lickorish \cite{Lic} proved that if we consider a link instead of a knot, we can get any 3-manifold from the 3-sphere in this way.
This tool has been heavily used in both knot theory and 3-manifold topology.

Thurston introduced hyperbolic geometry into the theory of Dehn surgery.
First of all, his uniformization theorem for Haken manifolds implies that any knot that is neither a satellite knot nor a torus knot has a complement which has a complete hyperbolic metric.
Such knots are called hyperbolic knots.
For a hyperbolic knot $K$, Thurston considered Dehn surgeries along $K$, and proved that except for  finitely many slopes, the manifolds obtained by surgeries are all hyperbolic.
(This theorem is called the hyperbolic Dehn surgery theorem.)\index{hyperbolic Dehn surgery theorem}
\index{hyperbolic Dehn surgery}
An interesting feature in the proof of this result is that it does not use the uniformization theorem, once we know the complement of $K$ has a complete hyperbolic metric.
Also among those manifolds obtained by hyperbolic Dehn surgery, there are non-Haken manifolds, whose hyperbolicity cannot be shown by the uniformization theorem.

More generally, by using the same techniques as those used in the proof of the hyperbolic Dehn surgery theorem, Thurston showed that for any complete hyperbolic 3-manifold of finite volume having a torus cusp, one can obtain a hyperbolic 3-manifold by attaching a solid torus to a cusp. Such an operation gives a hyperbolic 3-manifold except for finitely many homotopy classes of the attaching disc.
Combining this theorem with his result on the Gromov invariant, which we shall present in \S \ref{s:norm}, he obtained the fact that the set of the volumes of hyperbolic 3-manifolds constitues an ordered subset of $\mathbb R$ isomorphic to $\omega^\omega$.
The volumes\index{volume} of hyperbolic 3-manifolds are important objects in 3-manifold topology and there is still  a large amount of activity taking place on this topic.

Hyperbolic Dehn surgery\index{hyperbolic Dehn surgery} gives a better framework to understand deformations of hyperbolic cone structures.
The homotopy classes of simple closed curves on a torus can be regarded as co-prime lattice points on $\mathbb R^2$.
Therefore, the set of Dehn surgeries on a hyperbolic knot $K$ can be identified with such lattice points.
In this picture hyperbolic cone structures on $S^3$ whose singularities are isotopic to $K$ can be identified with a subset of the $x$-axis.
In this respect, Thurston considered a hyperbolic Dehn surgery space which contains both hyperbolic cone structures and hyperbolic Dehn surgeries.

                            \subsection{Kleinian groups}
                            \label{Kleinian}
                            
The notion of Kleinian\index{Kleinian group} group was first introduced by Poincar\'{e} as a generalization of the notion of Fuchsian group.\index{Fuchsian group}
Kleinian groups were extensively studied from the viewpoint of complex analysis in 1960-1975 by Ahlfors, Bers, Kra, Maskit, Marden among others.
In particular, Bers considered the compactification of a slice in  quasi-Fuchisan space, which is now called the Bers compactification\index{Bers compactification} of Teichm\"{u}ller space \cite{Bers}, and Marden considered the deformation space of convex cocompact representations \cite{Marden}.

In the process of proving the uniformization theorem for Haken manifolds,\index{uniformization theorem!Haken manifolds}  Thurston needed to develop the theory of deformation spaces of Kleinain groups\index{deformation spaces!Kleinain groups} and ends of hyperbolic 3-manifolds.\index{end of a hyperbolic 3-manifold}
%to show some sequences of Kleinian groups, which are images of the same group, to converge as representations, and analyse the property of limit Kleinian groups.
Thurston's proof of the uniformization theorem is largely divided into two cases, the one where the manifold is fibered over the circle and the other when it is not.
In the first case, treated in \cite{Th2}, he proved what is called the double limit theorem.\index{double limit theorem}
In the second case, he showed the compactness theorem for deformation spaces of hyperbolic structures on acylindrical 3-manifolds, proved in \cite{Th1}, and its relative version, proved in \cite{Th3}.
All these theorems are very important in the study of deformation spaces. They were generalized and became the fundamental tools in Kleinian group\index{Kleinian group} theory.
The reader may refer to Chapter 8 of the present volume,  by Lecuire \cite{Lecuire}. 

Another important part in Thurston's proof is the analysis of geometrically infinite Kleinian groups\index{geometrically infinite Kleinian group}\index{Kleinian group!geometrically infinite} in which
Thurston introduced the notion of geometric tameness using pleated surfaces.
A geometrically infinite (torsion-free and finitely generated) Kleinian group is said to be geometrically tame\index{Kleinian group!geometrically infinite!gemetrically tame} when every geometrically infinite end of the corresponding hyperbolic 3-manifold has a sequence of pleated surfaces tending to it.\index{geometrically tame}
Thurston showed that if a Kleinian group is geometrically tame, then the corresponding hyperbolic 3-manifold can be compactified by adding a boundary component to each of its ends.
This property is called the topological tameness.

Thurston considered that one can prove Marden's tameness conjecture\index{Marden's tameness conjecture} \cite{Marden} saying that every complete hyperbolic 3-manifold with finitely generated  fundamental group is topologically tame, by showing that every Kleinain group is geometrically tame.
He gave a proof of this conjecture in the special case of groups that are are algebraic limits of quasi-Fuchsian groups.

In \cite{ThB}, Thurston gave a list of 13 problems on Kleinian groups, one of which is Marden's tameness conjecture described above.
We shall say more about these problems in \S \ref{Kleinian impact}.

%It was by the work of J{\o}rgensen who proved the existence of hyperbolic orbifold structure on a surface bundle over the circe whose fiber is a torus with one cone point.

\subsection{The Thurston norm, the Gromov norm and the Gromov invariant}\label{s:norm}
As the work of Waldhansen, Haken, Jaco, Shalen and Johannson show, incompressible surfaces\index{incompressible surface} are important tools to study 3-manifolds.
In the same vein, understanding Seifert surfaces is essential for knot theory.
For every knot in the 3-sphere, the complement has first homology group isomorphic to $\mathbb Z$, and its Poinca\'{e} dual is represented by a Seifert surface.\index{Seifert surface}

Thurston introduced in \cite{T-Norm} a pseudo-norm on the second homology groups (or on the first cohomology groups) of 3-manifolds, which is called today the Thurston norm.\index{Thurston norm}
Given a second homology class $c$ of the manifold, its Thurston norm is defined to be $x(c)=\max\{\min_S( -\chi(S)), 0\}$, where $S$ ranges over all surfaces representing this homology class.
(In the case of 3-manifolds with boundary, it is more reasonable to consider their relative homology groups.) 
It should be noted that a surface realizing the Thurston norm is always incompressible, but the converse does not hold in general.

This simple idea led to very interesting results.
The norm can be extended first to the second homology with rational coefficient by defining $x(c/r)$ to be $x(c)/r$ for any integral homology class $c$, and by continuity to homology with real coefficients.
The unit ball in the second (relative) homology group with real coefficients constitutes a convex polytope with vertices on lattice points when the manifold is irreducible, atoroidal, and acylindrical.
This is called the Thurston norm polytope.\index{Thurston norm!polytope}

In a 3-manifold which fibers over the circle with fibers having negative Euler characteristic, the Euler class of the bundle of planes tangent to the fibers defines a second cohomology class.
Thurston proved that the Thurston dual norm\index{Thurston norm!dual} of such an Euler class is always equal to $1$.
This implies that for any fibering over the circle, the second homology class represented by a fiber lies in the interior of a facet of the Thurston norm polytope.
As a corollary, this implies that a 3-manifold fibered over the circle (whose fibers have negative Euler characteristic) contains an incompressible surface which cannot be a fiber for any fibering over the circle.
More generally, if a 3-manifold has a transversely oriented codimension-one foliation, we can consider the Euler class of the bundle of planes tangent to the leaves.
Thurston showed that if this foliation does not have Reeb components, then the Thurston dual norm of the Euler class is less than or equal to $1$, and that in particular if such a foliation has a compact leaf, then the dual norm is equal to $1$.
He also showed that a compact leaf of such a foliation realizes the Thurston norm of its second homology class.
The notion of Thurston norm and these results gave the foundation of further studies of fibrations and foliations in 3-manifolds.

We now consider the Gromov norm.
             
 To introduce the notion of bounded cohomology in \cite{Gr}, Gromov defined a pseudo-norm on homology groups as follows.\index{Gromov norm}
 Given a singular chain $s = \sum a_i c_i$ with real coefficients, define  its norm $\Vert s \Vert$ to be $\sum_i |a_i|$.
  Then for any homology class $\sigma$, its Gromov norm is defined to be $\inf \{\Vert s \Vert \mid [s]=\sigma\}$.
 The same definition works also for relative homology groups.
  In particular, for a closed manifold, the Gromov norm of its fundamental class is called its Gromov invariant.\index{Gromov invariant}
  Gromov proved that for a closed manifold of dimension $n$, its Gromov invariant is equal to its volume divided by a constant $v_n$ depending only on $n$, where $v_n$ is equal to the supremum of the volumes of $n$-dimensional hyperbolic simplices.

  Thurston took up Gromov's invariant as the topic of Section 6 of his lecture notes \cite{T-Notes3}.
  (This was before the publication of Gromov's paper \cite{Gr}, in which we can find Thurston's influence, both explicitly and implicitly.)
 Thurston generalized Gromov's result to negatively curved manifolds, where the Gromov norm is bounded below by the volume divided by a constant depending only on $n$, and to manifolds with geometric structures, where equality between the Gromov norm\index{Gromov norm} and volume holds, but with a constant depending on the geometry.
 Thurston also considered a relative version of the Gromov invariant for manifolds with boundary.
 Using this in the case of hyperbolic 3-manifolds with torus cusps, he showed that the operation of hyperbolic Dehn surgery decreases the volume.
 He also proved that a torus cusp corresponds to an accumulation point of the distribution of volume, and that this is the only way for volumes of hyperbolic 3-manifolds to accumulate.
 
 The Gromov norm for the second homology group of a 3-manifold is related to the Thurston norm.
 Thurston conjectured that if we change embedded surfaces in  the definition of the Thurston norm to immersed surfaces, then the norm  obtained should coincide with the Gromov norm.
 Gabai proved in \cite{Gabai1982} that this is indeed the case.

  \subsection{Conformal geometry and holomorphic dynamics}\label{s:conformal}

Conformal geometry,\index{conformal geometry}\index{holomorphic dynamics} since its birth,  is intertwined with topology. The relation between the two subjects started in Riemann's doctoral dissertation  (1851) \cite{Riemann} in which  he introduced the concept of Riemann surface (as we call it today), as a branched 
coverings of the complex plane or of the Riemann sphere. His work was partly motivated by the problem of describing a general method for finding domains of definitions for multivalued functions $w(z)$ of a complex variable $z$  defined by algebraic equations of the form $f(w,z)=0$, so that a multi-valued function becomes single-valued. (This is the original meaning of the word ``uniformization".)\index{uniformization} Thanks to the work of Riemann, analytic functions became objects that are no more necessarily defined on the complex plane or on subsets of it.  With him, the concept of Riemann surface,  with the closely related notion of analytic continuation, were born. Riemann further developed his
ideas on this topic in his paper on Abelian functions \cite{Riemann-Abelian}.    
  At the same time, he introduced  a number of topological notions that he used in the theory of functions of a complex variable: connectedness,
degree of connectivity,  genus, the classification of closed orientable surfaces, etc.
 One should also remember that there was still no notion of manifold in those times, and a surface could not be simply defined as a 2-manifold. 
 Riemann also studied moduli of Riemann surfaces. He discovered that the number of such moduli, for a  surface of genus $g$, is $3g-3$.  He also proved the famous \emph{Riemann mapping theorem},\index{Riemann mapping theorem} saying that any simply connected open subset of the complex plane, provided it is not the whole plane,  is biholomorphically equivalent to the open unit disc.  The problem of characterizing topologically analytic functions arose (this was called later on the ``Brouwer problem").  This problem was also a motivation for the development of quasiconformal mappings. Indeed, Gr\"otzsch, Lavrentieff and Teichm\"uller,  the three founders of the theory, tried to prove for quasiconformal mappings\index{quasiconformal mapping} some results that were known to hold for  conformal mappings, with the idea that it was only the topological form of the mapping that matters and not the fact that it was conformal.
 
 It was natural that Thurston got attracted by this field. 
   One of the first problems that he asked, when he joined the community of MathOverflow is related to Riemann surfaces and rational functions on the Riemann sphere. He wrote the following (posted on September 10, 2010): ``I would like to understand and compute the shapes of rational functions, or equivalently, ratios of two polynomials, up to Moebius transformations in both domain and range." He also  formulated the following more precise problem:  
``Given a set of points to be the critical values [in the range], along with a covering space of the complement homeomorphic to a punctured sphere, the uniformization theorem says this Riemann surface can be parametrized by $S^2$, thereby defining a rational function. Is there a reasonable way to compute such a rational map?"

On holomorphic dynamics, we shall mention two results of Thurston. The first one has a  discrete character; it is Thurston's topological characterization of rational maps among branched coverings of the sphere. The second one has a continuous character; it is his result with Sullivan on holomorphic motions.

  Thurston's theorem on the characterization of rational maps, that he proved at the beginning of the 1980s was a preliminary (but huge) step towards the program he formulated 30 years later on MathOverflow. It  gives a necessary and sufficient condition for a branched covering of the sphere which is  postcritically finite,\index{postcritically finite map} that is, such that the union of the forward orbits of the critical points is finite,    to be homotopy equivalent to a rational map.  Here, homotopy equivalence is defined in an appropriate and natural sense; the relation is called now Thurston's equivalence. Thurston's criterion is given in terms of the action of the  covering map (by taking inverse images) on systems of homotopy classes of essential simple closed curves on the sphere with the postcitical set deleted. The necessary and sufficient condition refers to this action, and it now carries the name ``absence of a Thurston obstruction." 
  
  Thurston obtained this theorem in 1982. He lectured on it on several occasions and he wrote notes that were widely circulated \cite{T-Notes-rational}.  Although it was announced at some point that a final version of these notes will be published in the CBMS conference series of the AMS, the notes remained unpublished. Adrien Douady and John Hubbard wrote a proof of Thurston's theorem, following his outline, and they circulated it in  preprint form in 1984. Their paper was eventually published in 1993 \cite{DH}.  
  
    The proof of Thurston's theorem, like the proofs of some others of his geometrical important results, involves the construction of a weakly contracting self-map of a certain Teichm\"uller space\index{Teichm\"uller space} (which in the case at hand is the one of the sphere with the post-critical set deleted). The main step  in this proof is to show that in the absence of a Thurston obstruction, this map has a unique fixed point. The fixed point, when it exists, is the desired rational function. 
    
    There are several analogies between this theory of Thurston and his classification theory of surface homeomorphisms: the use of hyperbolic geometry, the construction of an action on Teichm\"uller space, the study of an action on the collection of homotopy classes of simple closed curves, the existence of invariant laminations, the use of quasiconformal mappings, the  utilization of linear algebra in the study of the action on curves, and in particular the Perron--Frobenius theorem for nonnegative matrices, etc.

    Many applications of Thurston's theorem were obtained by various authors. We mention as an example its use in the theory of mating of two polynomials (getting a rational map whose Julia set is obtained by gluing those of two postcritically finite polynomials).
   We refer the reader to the paper \cite{BCT} by Xavier Buff, Guizhen Cui, and Lei Tan for a survey of Thurston's theorem, including a self-contained proof of a slightly generalized version of this theorem and an overview of its applications.   
         
As another aspect of holomorphic dynamics that was touched upon by Thurston, we review now his result with Sullivan on holomorphic motions. 
 
 A holomorphic motion of a subset $X$ of the complex plane $\mathbb{C}$  is a family of mappings $f_t:X\to \mathbb{C}$ parametrized by a complex number  $t$ (considered as a complex time parameter) varying in a domain $T$ containing the origin and satisfying the following three properties: (1) for each fixed $x$, $f_t(x)$ is holomorphic in $t$; (2) for each fixed $t$, $f_t(x)$ is injective in $x$; (3) $f_0$ is the identity mapping of $X$.  The motivation behind this definition is the wish to adapt the topological notion of isotopy to a holomorphic context. The main question addressed is to know whether a holomorphic motion of the subset $X$ extends to a holomorphic motion of the complex plane $\mathbb{C}$. (This is a holomorphic analogue of the topological problem of extending an isotopy to an ambient isotopy.)      Holomorphic motions\index{holomorphic motion} were introduced by Ma\~{n}\'e, Sad and Sullivan in \cite{MSS}.

  The main result of the paper \cite{ST1} by Sullivan and Thurston  says that there exists a universal constant $a>0$ such that any holomorphic motion of any subset $X$ of $\mathbb{C}$ parameterized by the unit disk $\{$t$:|t|<1\}$ can be extended to a holomorphic motion of the complex plane with time parameter the disc $\{$t$:|t|<a\}$.  
With his logician bias, Thurston did not fail to notice a close relation between the problem of extending holomorphic motions and a ``holomorphic axiom of choice".  This result was later improved by Slodowski \cite{Slodowski}, who showed in particular that one may take $a=1$ in the above statement, answering a question raised by Sullivan and Thurston in their paper, and proving a conjecture they formulated, precisely, related to the holomorphic axiom of choice.

 In the same paper, Thurston and Sullivan introduced the notion of quasiconformal motion.\index{quasiconformal motion} They noticed that the map $f_t$ in the above definition is necessarily quasiconformal and extends to a quasiconformal map of $\mathbb{C}$. They also proved a general extension theorem for quasiconformal motions over an interval. 
       In proving their results,  Sullivan and Thurston  introduced an averaging procedure for pairs of probability measures defined on the Riemann sphere.

  It was known, since the work of Gr\"otzsch, Lavrentieff and Teichm\"uller that quasiconformal mappings are useful in the study of conformal mappings. From the work of Sullivan and Thurston, the notion of quasiconformal motion became useful in the study of holomorphic motion. 
  
   There are applications of holomorphic motion to the theory of Kleinian groups, where the subset $X$ in the above definition of holomorphic motion is the limit set of the group action. There are also applications in  the theory of iterations of rational maps (where $X$ is taken to be the Julia set of the map), in the theory  of invariant metrics in complex geometry, in the study of holomorphic families of Riemann surfaces, in the theory of quasiconformal mappings and in the study of Teichm\"uller spaces.

    \subsection{Complex projective geometry} \label{s:grafting}
    
Complex projective geometry is a classical topic, rooted in  the 19th-century work of Klein, Poincar\'e and their contemporaries. One must add   that in the  1960s, Bers and his collaborators were thoroughly involved in the relation between the complex projective geometry of surfaces and  Teichm\"uller theory.  In particular, the Bers embedding of Teichm\"uller space is defined in the setting  of complex projective structures.

 In the late 1970s, Thurston reconsidered this theory. He did not publish any paper on this topic, but many authors wrote on it, following Thurston's outline.
He highlighted a profound analogy between the complex projective geometry of surfaces\index{complex projective surface} and the theory of Kleinian groups, opening a new perspective in 3-manifold topology, and motivating  important later works by  Sullivan, Epstein, Marden and others.
        
 Thurston introduced  metrics that are conformal to the complex structures on complex projective surfaces.\index{Thurston metric!projective surface} 
 In the simplest case, such a metric is obtained by grafting a Euclidean annulus on a hyperbolic surface, after cutting it along a simple closed geodesic. (In fact, it is the projective structure underlying the annulus, and not its Euclidean structure, which matters.) The general definition needs Thurston's theory of projective laminations, and the  metric is obtained as a limit of metrics when a sequence of simple closed curves along which the grafting is made converges, in Thurston's topology, to a measured lamination.
 This metric is now called Thurston's metric on the complex projective surface.\index{Thurston metric}  Thurston also gave a description of this metric as a Finsler metric, in fact, as a solution of an extremal problem,  in the spirit of the Kobayashi metric on a complex space: the distance between two points is the infimum of the length of piecewise $C^1$ paths joining them. Here, the length of a $C^1$ path is computed as the integral of the norms of vectors tangent to this path, the norm of a   vector $v$ being the infimum of the norms of vectors $v'$ in the tangent space of  the unit disc equipped with the hyperbolic structure, such that there exists a projective map from the unit disc to the surface whose differential sends $v'$ to $v$. 
 
 Before Thurston came into the subject, grafting, in its simplest form, was studied by Bers \cite{Bers}, Maskit \cite{Maskit} and Hejhal \cite{Hejhal}. There are also relations with the theory of harmonic maps between surfaces, in particular, between the extremal length of a grafted surface and the energy of a harmonic map, see \cite{Tanigawa}. The grafting operation makes connections between projective structures and hyperbolic geometry. Such connections were already known to Cayley, Klein, Study and others, and traces of the elementary grafting operation can be found in the work of Klein. 
 
 Using  the notion generalized grafting, Thurston discovered a geometric parametrization of the moduli space of marked projective structures on a surface,  establishing a precise relationship between this moduli space and the Teichm\"uller space of the surface. The parametrization takes the form of a homeomorphism between the moduli space of projective structures and  the product of measured lamination space\index{measured lamination space} with Teichm\"uller space. The result says that the fiber at each point of the natural forgetful map from the moduli space of projective structures to Teichm\"uller space is naturally identified with the space of measured laminations on the surface. 

 There are several surveys of Thurston's work on complex projective geometry, and they give complementary points of view on the subject.  The paper \cite{KaT} by Kamishima and Tan is concerned with grafting and Thurston's parametrization in the setting of the theory of geometric structures and locally homogeneous spaces. Goldman, in his paper \cite{Goldman-P}, sets the foundations of the theory of complex projective structures on surfaces as geometric structures, using the notions of holonomy and developing map, in the tradition of Thurston, and with ideas originating in the work of Ehresmann \cite{Ehresmann}. He refers, for Thurston's parametrization of the moduli space of projective structures by the product of measured lamination space with Teichm\"uller space, to a course given by Thurston at Princeton University in 1976--1977 (Goldman was an undergraduate student there). 
 We also refer the reader to the paper \cite{Dumas2007} by Dumas, and to his survey \cite{Dumas-H}, for the parametrization of the moduli space of complex projective structures on surfaces.   Chapter 6 in the present volume, by Shinpei Baba \cite{Baba}, is a survey of Thurston's work on complex projective structures on surfaces, and it contains other references.
 
 We also mention the paper \cite{Sullivan-T-canonical} by Sullivan and Thurston, in which these authors provide a series of examples that show the subtleties of higher-dimensional real and complex projective structures, together with other kinds of geometric structures (inversive and affine). 

Now we wish to talk about a classical parametrization of the moduli space of complex projective structures on surfaces based on the Schwarzian derivative, a differential operator invariant under M\"obius transformations.  This parametrization first appeared in the nineteenth century in the works of Schwarz, Klein and others.  Thurston used it in his study of projective geometry, and we briefly discuss this.

In his paper \cite{T-zippers} published in a  special volume of the AMS at the occasion of the proof of the Bieberbach conjecture, Thurston introduced a topology on the set of  univalent mappings from the unit disk into the Riemann sphere using the topology of uniform convergence of Schwarzian derivatives. The uniformity refers to the hyperbolic metric of the disk. To see that this is a natural object of study, one may recall that the Schwarzian derivative\index{Schwarzian derivative} was classically used to obtain Riemann mappings of some special domains of the complex plane: regular polygons, domains with circular edges, etc.; generalizing, this makes the set of Schwarzian derivatives is a candidate for a natural parameter space for projective structures.

 The usual definition of the Schwarzian derivative, involving  third-order complex derivatives, makes it at first sight quite obscure. It is interesting to see how this object was described by Thurston in his paper. He writes: 
``For the benefit of people to whom the Schwarzian derivative may seem a mystery, we will set the stage by discussing the Schwarzian derivative." He continues a few lines below:  ``The Schwarzian derivative is very much like a kind of curvature: the various kinds of curvature in differential geometry measure deviation of curves or manifolds from being flat, while the Schwarzian derivative measures the deviation of a conformal map from being a M{\"o}bius transformation." Then, after a page of explanations, he writes: ``A formula for the Schwarzian derivative can be readily determined from the information above, or it may be looked up---someplace else. Like the formula for the curvature of a curve in the plane, the formula looks somewhat mystical at first, and in a quantitative discussion the formula tends to be a distraction from the real issue." Reading Thurston's paper gives a clear intuition of what the Schwarzian derivative is.

Responding to a question by Paul Siegel on MathOverflow, on September 9, 2010, asking: ``Is there an underlying explanation for the magical powers of the Schwarzian derivative?", Thurston writes, on the next day: ``Like many people (but not all people), I have trouble thinking in terms of formulas such as that for the Schwarzian. For me, a geometric image works much better [\ldots]". He then gives a nice description based on hyperbolic geometry and quadratic differentials and measured foliations.
On the next day (September 11, 2010), Thurston writes  to   Siegel who thanked him for his response: ``I appreciated the question, which resonated with my thoughts. I'm new to MO, but it seems like a rich environment. I understand MO is not intended for extended threads, but I'd like to leave a pointer forward to my first question, which I posted partly as a followup to this, since it indicates the immediate source for my interest in Schwarzians." In this follow up, Thurston asks several questions, including the following: 
\begin{quote}\small
Given a set of $2d-2$ points on $\mathbb{C}P^1$ to be critical points [in the domain], it has been known since Schubert that there are Catalan rational functions with those critical points. 

--- Is there a conceptual way to describe and identify them?

--- Is there a complete characterization of the Schwarzian derivative for a rational map, starting with the generic case of $2d-2$ distinct critical points?

--- What planar graphs occur for Schwarzian derivatives of rational functions? What convex (or other) inequalities do they satisfy?
\end{quote}

                             \subsection{Circle packings and discrete conformal geometry}\label{s:circle}

   The study of circle packings, that is, configurations of circles that are tangent to each other, is classical and can be traced back to the work of Apollonius of Perga (3d century B. C.),  see \cite{Apollonius}.\footnote{For  Apollonius' works, the main reference is  Rashed's critical edition of the Arabic manuscripts (many Greek textes do not survive), published by de Gruyter in 5 volumes (more than 2500 pages) between 2008 and 2010. Apollonius' problems\index{Apollonius problems} are discussed in the volume \cite{Apollonius}.} 
 In the nineteenth century, Paul Koebe proved the existence of some circle packings, and  considered  the idea of using them to prove the Riemann Mapping Theorem \cite{Koebe}. For a review of the work of Koebe, we refer the reader to the chapter by Philip Bowers  in the present volume \cite{Bowers}. One may recall in passing that Koebe (and independently Poincar\'e) proved a wide generalization of the Riemann Mapping Theorem, namely, the uniformization theorem.\index{uniformization theorem}

  Thurston's work on circle packings inaugurated the notion of discrete Riemann mapping theorem,\index{Riemann mapping theorem!discrete} and more generally, the study of discrete conformal mappings.\index{discrete conformal mapping} At the same time, it shed a new light on several  geometric ideas that are rooted in classical mathematics.

 Let us first recall that the (classical) Riemann mapping theorem,\index{Riemann mapping theorem} proved by Riemann in his doctoral dissertation \cite{Riemann}, says that an arbitrary simply connected open subset of the complex plane which is not  the whole plane  is conformally equivalent to the unit disc by a mapping which is unique up to composition by a M\"obius transformation.

Thurston conjectured the existence of a discrete version of the Riemann mapping theorem\index{Riemann mapping theorem!discrete} as a limit of a sequence of circle packings. The intuition behind this is that a conformal mapping between two open subsets of the plane is characterized by the fact that it sends infinitesimal circles to infinitesimal circles (recall that at the level of tangent planes,  it sends circles centered at the origin to circles centered at the origin). Therefore one might hope that finding circle packings with smaller and smaller radii on a given domain and a sequence of homeomorphisms that send them to circle packings of the unit disc  leads, by taking limits, to a  Riemann mapping. A ``discrete Riemann mapping"\index{discrete Riemann mapping theorem}\index{Riemann mapping theorem!discrete} is one of these mappings used in the approximating sequence.

                   Thurston's conjecture was proved by Rodin and Sullivan in their paper \emph{The convergence of circle packings to the Riemann mapping}  \cite{RoS}.  In their introduction, the authors recall the setting: \begin{quote}\small
                   In his address,\footnote{International symposium in Celebration of the Proof of the Bieberbach Conjecture. Purdue University, March 1985.} \emph{The finite Riemann Mapping Theorem}, Bill Thurston discussed his elementary approach to Andreev's theorem and gave a provocative, constructive, geometric approach to the Riemann mapping theorem. This method is quite beautiful and easy to implement on a computer.
                   \end{quote}
                    They then recall Thurston's strategy of the proof:
                    \begin{quote}\small
                    Almost fill a simply connected region $R$ with small circles from a regular hexagonal circle packing. Surround these circles by a Jordan curve. Use Andreev's theorem to produce a combinatorially equivalent packing of the unit disc---the unit circle corresponding to the Jordan curve. The correspondence between the circles of the two packings ought to approximate the Riemann mapping.
                    \end{quote}
                                     
                                      Following Thurston's ideas, Rodin and Sullivan develop in an appendix to their paper an algorithm to obtain the discrete Riemann mapping.
                                      
                                In his Princeton lectures, Thurston studied circle packings in the midst of a discussion of orbifolds and of an existence theorem for hyperbolic polyhedra. 
When he started lecturing on the subject, he was  aware neither of Koebe's nor of Andreev's work; see the interesting historical remarks in Bowers' review \cite{Bowers-review}.  He realized at some point that some of the results he obtained were  generalizations of results contained in two papers by Andreev \cite{Andreev1, Andreev2}.  He then called the existence theorem for circle packings that is contained in Chapter 13 of his Princeton notes \cite{T-Notes3}  Andreev's Theorem. The result is now called the\index{Koebe--Andreev--Thurston theorem} Koebe--Andreev--Thurston theorem.  This theorem states that for any triangulation of a closed surface (of any genus) which lifts to a simple triangulation of the universal cover (that is, a triangulation which has no pair of edges connecting the same vertices, and no edge connecting the same vertex), there  exists a unique metric of constant curvature on the surface with a circle packing that is modeled on it. Furthermore, the packing is unique up to a conformal map isotopic to the identity (which implies, in the hyperbolic case, that the map is the identity). Thurston deduced the uniqueness result from  Mostow's rigidity theorem.
 In his notes, he considered in detail the genus 1 and $\geq 2$ cases. The genus $0$ case was treated by Andreev and was attributed to him by Thurston. Marden and Rodin wrote a paper showing that Thurston's method also gives a proof for the genus 0 case \cite{MR1990}.
 
  Thurston also proved an existence theorem for patterns of circles that generalizes a result of Koebe, allowing an overlap among the pattern of circles (and not only tangency), and he used this result in his proof of the generalized  Andreev theorem.  In  \S 13.4 of  his notes \cite{T-Notes3}, he gave algorithms for constructing circle packings and more general circle patterns. His algorithms allow computations.

      For an overview of Thurston's discrete Riemann mapping theorem and its impact, we refer  the reader to the comprehensive survey by Bowers in the present volume \cite{Bowers}.
      We also refer to Luo's paper \cite{Luo-R} and to Kojima's survey  \cite{Kojima}.

       \subsection{Word processing in groups, automata and tilings}
       
 Besides the name of Thurston, two names will be highlighted in this section: Jim Cannon and John Conway.

We start with groups and automata, to which the name of Cannon is attached.     
      
 In 1984, Cannon published a paper in which he showed that Cayley graphs of cocompact discrete groups of isometries of $n$-dimensional hyperbolic space can be given finite recursive descriptions  \cite{Cannon}.  He wrote in the introduction that he was inspired by Thurston, who showed that a large number of groups that are of interest to topologists cannot be dealt with using the standard methods of combinatorial group theory, but can be attacked by ``a return to geometric consideration", that is, the classical methods of Dehn and Cayley.  It is significant that Cannon's paper contains an appendix on elementary properties of hyperbolic space, for which, at that time, no modern exposition was available, except for Thurston's unpiblished notes \cite{T-Notes3}.

    Thurston noticed that Cannon's results can be formulated in the language of finite state automata, and may be applied to a wider class of groups. This led him to the introduction of the notion of automatic group.\index{automatic group} This is a group equipped with a simple algorithm for the word problem, that is, an automaton can tell when two words (in a given system of generators) represent the same element in this group.   
        
    After
 their discovery by Thurston, automatic groups found applications in a wide class of domains including low-dimensional topology and geometry, geometric and combinatorial group theory,    algorithmics, decision theory,  computer vision, mathematical logic, etc. Furthermore, the theory of automatic groups is closely related to that of finite state automata, which has applications in computer science. Thurston was interested in all these applications. He developed with his collaborators computer programs to carry out what he called ``word processing on groups." There is also a relation with self-similar tilings. Thurston writes in \cite{T-AMS1989}: ``An automatic structure for a group in general produces a kind of self-similar tiling of a certain `sphere at infinity' for the group; in particular examples, this space is actually a 2-sphere."    

      Soon after their discovery, automatic groups became a central part of geometric group theory.
Thurston collaborated with Cannon and others on this theory in connection to his eight geometries of 3-manifolds. In the paper \cite{CGFT} written with Cannon, Floyd end Grayson, he showed that no cocompact discrete group based on solvgometry, {\bf Sol}, is almost convex. As a consequence, the Cayley graph of such a group cannot be efficiently constructed by finitely  local replacement rules. After recalling Thurston's geometrization conjecture,\index{geometrization conjecture} the authors write that ``any package of decision algorithms designed to compute within the fundamental groups of low-dimensional manifolds and orbifolds must be able to deal with the groups from each of the standard geometries." The theory of automatic groups is developed in the comprehensive textbook \cite{T-Automatic}  that Thurston wrote with Epstein, Cannon, Holt,  Levy and Paterson, which appeared several years after he started working on this topic.  

We now pass to tilings.\index{tiling} 

The  study of tilings is closely related to discrete group actions, a theory  that plays an essential role in Thurston's work on 3-manifolds. Thurston was fascinated not only by tilings in dimension 3, but also by the theory of plane tilings, in particular by the theory of quasiperiodic tilings. These are tilings where finite patters appear regularly, without being necessarily periodic. 
Let us quote an excerpt   from a set of questions that Thurston distributed at the beginning of his course on ``Geometric topology" at Princeton, during the Spring Semester of 1983 \cite[Question 9]{T-Questions-1983}: 

 \begin{quote}\small
 Is there a general mathematical theory for Penrose-like tilings,\index{Penrose-like tiling} where one specifies certain combinatorial relationship and then deduces that certain shapes of tiles exist which satisfy these relations? Are there many essentially different such tilings, or just few?
 \end{quote}

   Thurston was stimulated on this subject by ideas of  Conway, who was working at the same university and who made major contributions to group theory, sphere packings, tilings and cellular automata.  
 The latter, together with Jeffrey Lagarias, developed a method, based on combinatorial group theory,  to tackle the problem of tiling some finite region of the plane using a certain number of regular tiles \cite{Conway-L}. 
This method involved the encoding of the edges of the tiles by elements of a finitely-presented group in such a way that a tile can be interpreted as a relator in the group. The problem was then reduced to that of deciding whether some group element, describing the boundary of the plane region, is the trivial element.  

In 1990, Thurston published a paper in the American Mathematical Monthly \cite{T-Conway} in which he re-interpreted Conway's construction using the tools of geometric group theory.  In the same paper, he gave a necessary and sufficient condition for a simply-connected region of the plane which is the union of unit squares, to have a tiling by dominos, that is, rectangles which are the union of two squares. 
 He also gave several constructions of tessellations of planar regions by given tiles.
 
 In the same year, Thurston and Conway, together with Peter Doyle, started a new course at Princeton, called ``Geometry and the imagination." Thurston writes  \cite{T-Math-edu}: ``The course came alive, qualitatively more than any course we had taught before.
Students learned a lot of mathematics and solved problems we wouldn't have dared
ask in a conventional college class."

Thurston's collaboration with Conway includes the paper  \cite{CDHT}  by Conway,  Delgado Friedrichs, Huson and Thurston in which these authors gave  a new enumeration of $n$-dimensional crystallographic 
 groups, that is,  cocompact discrete subgroup of the isometry group of Euclidean 3-space. The enumeration is based on a description of these groups as fibrations over the plane 
crystallographic groups, when the enumeration is possible. The ``old" enumeration, due to Barlow, Fedorov, and Sch\"{o}nflies, dates back to the 1890s.

We mention  another paper on tilings (although this word is used in a slightly different meaning), namely, the  paper \cite{CTST}, by Coven,  Geller, Silberger and Thurston, concerned with the symbolic dynamics of tiling the integers. Here, a finite collection of finite sets of integers is said to ``tile the integers" if  the set of all integers can be written as a disjoint union of translates of elements of this finite set.    These elements are called tiles. To such a set of tiles, the authors associate a bi-infinite sequence of elements of tiles. They show that the set of all such sequences  is a sofic system, and that every shift of finite type can be realized (up to a power)  as a tiling system.
 
 The paper \cite{T-AMS1989} contains results on self-similar tilings,\index{self-similar tiling} in particular, constructions of such tilings from algebraic integers $\lambda$ whose Galois conjugates, except $\lambda$ and $\overline{\lambda}$, are smaller. More generally, Thurston introduced the notion of  complex expansion constants for self-similar tilings, and he gave a characterization of these constants. He obtained a characterization of the set of similarities for self-similar tilings of the plane or of higher-dimensional spaces, making an analogy with the construction of Markov partitions from classical dynamical systems. Beyond the results he obtained, Thurston emphasized the aesthetical side of the topic. He writes: ``What is interesting about this subject is the particular constructions---at issue is how simple and how nice can a self-similar tiling can be."
 
  Automata\index{automaton} and tilings\index{tiling} were part of the subject of a series of lectures which Thurston gave at a summer AMS colloquium in the summer of 1989. The title of the series was \emph{Groups, tilings and finite state automata}. A preprint carrying the same title  \cite{T-AMS1989} was distributed at the meeting, and it was later included in the Research Reports of the Geometry Center preprint of the University of Minnesota, a center co-founded by Thurston.  The paper, which may be considered as semi-expository, remained in a preprint form. In this domain, and like many other topics which he considered, Thurston had a huge amount of ideas bubbling in his brain, and it was certainly difficult for him to sort out what was new and what was known in some sort or another.

 \subsection{Computers}
In the preceding section, we were led to mentioning computers\index{computers} quite a few times.\index{computer science} We give here a quick overview on other works of Thurston on this subject, and of his collaboration with computer scientists.
We highlight the fact that Thurston's collaboration with computer scientists was two-fold. On the one hand, he used methods of geometry, in particular 3-dimensional hyperbolic geometry, to solve problems in computer science, 
and on the other hand, motivated by questions that arose from computer science, he developed new topics and opened up new ways of research in geometry,.

Thurston was heavily involved in computing and computer graphics\index{computer graphics} since the 1980s. Let us quote a question from his list addressed to his students that we already mentioned \cite[Question 18]{T-Questions-1983}:

 \begin{quote}\small
 What is the information content of text? How well can one model the sequence of letters in a novel as a dynamical system? That is, suppose you forget that you know anything about language and meaning, and just try to analyse it from a statistical point of view; how could you do in automatically guessing what the next letter would be?
 
  This relates to the question of how much space it takes to store such a string of text in a computer.
  Given a model process, one could make up a coding scheme. In one direction, it would be possible to feed in a random set of bits and have the code produce a more-or-less plausible stream of text (depending on the complexity of the process which one allows); and in the other direction, one would feed in a text and have it compressed into a much shorter stream of symbols. One would try not to be prejudiced too much by the meaning of the  words, but still use knowledge of English (or whatever language) to figure out a good reasonably small set of data which are useful in predicting what next occurs.

  Similar questions can be asked about many other human-generated processes (e.g., music), many of them with obvious applications (e.g., the stock market, sequences of answers to multiple choice texts, \ldots). How much entropy do these processes have? Are there families of dynamical systems which do well the modeling? 

         \end{quote}

 Thurston collaborated with computer scientists on geometric problems he formulated, but also on problems that were asked by computer scientists themselves. We mention first his paper with Sleator and Tarjan, \emph{Rotation distance, triangulations, and hyperbolic geometry} \cite{Rotation}, published in 1988, in which a distance, called \emph{rotation distance},\index{rotation distance} is defined on the set of binary trees, as the minimal number of rotations that may be used to convert one of these trees into the other. The term ``rotation" denotes here the operation of collapsing an internal edge of a binary tree to a point and expanding the node, obtaining a new binary tree. The authors show in this paper that for binary trees\index{binary tree} with $n$ nodes with $n\geq 11$, the maximum rotation distance is at most $2n-6$. The motivation for this problem comes from a problem used in data structuring and network algorithms, and more precisely, from a conjecture called the \emph{splaying conjecture}.\index{splaying conjecture} The authors explain this conjecture as follows: ``Splaying
is a heuristic for modifying the structure of a binary search tree in such a way
that repeatedly accessing and updating the information in the tree is efficient." The methods used in this paper are based on hyperbolic geometry, in the pure Thurston tradition. The rotation operation between binary trees is converted to an equivalent operation of  flipping a diagonal in a polygon then passing one dimension higher which permits the rotation distance problem to be reduced to a 3-dimensional  problem of dissecting hyperbolic polyhedra into tetrahedra.  The  volume of hyperbolic polyhedra appears in various forms as a fundamental object in this study. The last section of the paper contains open questions, asking in particular for more calculations of triangulations and volumes for polyhedra. A relation with the Gromov norm, as a measure of how efficiently a 
homology class in a hyperbolic manifold can be represented by simplices (Chapter 6 of Thurston's Princeton notes \cite{T-Notes3}) is also mentioned in this section on open problems. 
 
 There are other papers of Thurston on computer science and algorithmic problems in which Thurston's geometrical methods are used as en essential tool. We mention the three papers in collaboration with Gary Miller, Shang-Hua Teng and Stephen Vavasis \cite{Miller1997},  \emph{Automatic mesh partitioning} \cite{Miller1993}, 
  \emph{Separators for sphere-packings and nearest neighborhood graphs} and   \emph{Geometric separators for finite-element meshes} \cite{Miller1998},   
 and his paper with Bob Riley, \emph{The absence of efficient dual pairs of spanning trees in planar graphs} \cite{Riley}. 
 
 Finally, we mention Thurston's paper  \emph{Shapes of polyhedra and triangulations of the sphere} \cite{Thurston-shapes} motivated by the question of classifying the combinatorial classes of triangulations of the sphere with at most 6 triangles at a vertex, in which he was led to endow the 
 moduli space of polyhedra with $n$ vertices with given total angles less than $2\pi$ at each vertex  (that is, Euclidean cone metrics of nonnegative curvature) with a K\"ahler metric which is locally isometric to complex hyperbolic space $\mathbb{CH}^{n-3}$.  This paper had an enormous influence and several generalizations of the results were attempted by many authors.

Regarding his collaboration with computer scientists, Thurston writes his 1987 \emph{Notices} article \cite{T-Funding}: 
\begin{quote}
\small
Recently, through circumstances, I have spent time with computer scientists. I find myself talking and thinking about computer science problems, and analyzing them with modes of thought sometimes foreign to the culture of computer science. I enjoy this. My experience would be similar if I were to spend time with physicists, biologists, economists, chemists, engineers \ldots.
\end{quote}

One should emphasize the fact that Thurston since the 1970s has been constantly programing, computing, implementing lists of knots,  of 3-manifolds, of volumes of hyperbolic manifolds,  of tilings, etc.

    \subsection{Surfaces, mapping class groups and Teichm\"uller spaces} \label{s:surfaces}

    In 1975-76, Thurston gave a course at Princeton on the geometry and dynamics of  homeomorphisms of surfaces. He presented there a complete theory which had to have a major and everlasting  impact on low-dimensional topology and geometry. A major role in this theory was  given to Teichm\"uller space,\index{Teichm\"uller space} the space of isotopy classes of metrics of constant curvature -1 on a surface.  Thurston's  results included a compactification\index{Teichm\"uller space!compactification} of this space by the space of projective classes of measured foliations, the latter seen as a completion of the set of homotopy classes of simple closed curves on the surface. The results also included a natural action of the mapping class group on this compactified Teichm\"uller space, and the classification of mapping classes into periodic, reducible and pseudo-Anosov, obtained by analyzing the fixed point set of the action of a mapping class on this compactified space. 
    
     Copies of a set of notes on Thurston's course, taken by Bill Floyd and Mike Handel, were circulated, and in particular they arrived to Orsay where they gave rise to the famous seminar \emph{Travaux de Thurston sur les surfaces} which took place during the year 1976-77; see the paper \cite{EMS-Surveys} for the history of this seminar.   It appears that Thurston was already thinking about surfaces, and in particular how simple closed curves approach a foliation, at the time he was a PhD student in Berkeley, see Sullivan's account in \cite{EMS-Surveys}.

    A couple of years after the Orsay seminar, the book \cite{FLP} was written and became the major reference on Thurston's  theory on surfaces. In the meanwhile,  Thurston wrote a research announcement of his results, which he submitted without success to a few journals. The research announcement eventually appeared in  
    the Bulletin of the AMS, in 1988 \cite{Thurston-FLP}, 12 years after Thurston wrote it. The paper contains new bibliographical references and a new preface in which Thurston gives a few notes on the history of the manuscrit and of the subject. 

    Shortly after Thurston obtained his classification theorem for mapping classes, Lipman Bers gave a new proof of that theorem in a complex analysis setting, and using the techniques of quasiconformal mappings\index{quasiconformal mapping} \cite{Bers1978}. Bers's proof also uses the action of the mapping class group on Teichm\"uller space, but unlike Thurston's proof, it is based on an analysis of the translation length of an element of the mapping class group with respect to Teichm\"uller's metric. In fact, in Bers' classification, there are four sorts of mapping classes, according to whether the translation length is zero or positive, and in each case, according to whether this translation length (which is defined as an infimum), it is attained or not by a point in Teichm\"uller space. 
    
    The book \emph{Travaux de Thurston sur les surfaces} did not include Thurston's theory of geodesic laminations and train tracks, which turned out to be very efficient tools in low-dimensional topology. These notions were expanded on in the courses that Thurston gave the following years at Princeton, and they are included in his notes \cite{T-Notes3}.  Several books appeared on the subject, see e.g. the notes by Casson and Bleiler that arose from a course that Casson gave on Thurston's theory of surfaces at the University of Texas at Austin, \cite{Casson-Bleiler} and the book \cite{PH} by Penner and Harer on train tracks.

   One consequence of Thurston's work was the revival of nineteenth-century hyperbolic geometry, a subject which was almost forgotten.
 Thurston's notes  \cite{T-Notes3}, together with the books  \cite{FLP} and \cite{Casson-Bleiler}, were for some time the main references on this topic. (In particular,  \cite{FLP} contains all the hyperbolic trigonometry formulae that are useful in the theory). At the time Thurston started working on the subject, there were practically no modern treatments of the subject. Of course, one could   look into Lobachevsky's works, but this was very unlikely. The textbooks by Beardon, Ratcliffe, Anderson and others appeared several years later. 
    The classical books \emph{Elementary geometry in hyperbolic space} \cite{Fenchel} by Fenchel and \emph{Discontinuous groups of isometries in the hyperbolic plane} \cite{FN}
    by Fenchel and Nielsen, which existed in the form of lecture notes and had trouble in being published, appeared in 1989 and 2003, after Thurston's work made them famous. The so-called Fenchel--Nielsen parameters for hyperbolic surfaces, associated with geodesic pairs of pants decompositions, acquired all their strength after Thurston used them in his work. Works of Abikoff \cite{Abikoff1980}, Wolpert \cite{Wolpert1982} and others on this deformation appeared after Thurston revived the subject.

Pseudo-Anosov homeomorphisms, which appeared in Thurston's classification, turned out to be a major ingredient in the geometry and topology of 3-manifolds; we mention for instance their role in the construction of hyperbolic manifolds which fiber over the circle, explained in \ref{3-man}.

              Before continuing on Thurston's theory of surfaces, we make a small digression concerning Nielsen's contribution to this subject.

   Jakob Nielsen, in several long papers published between 1927 and 1944 \cite{Nielsen1, Nielsen2, Nielsen3, Nielsen4}, studied  questions related to automorphisms of surfaces, using   hyperbolic geometry.
    Thurston writes in the introduction of his paper \cite{Thurston-FLP}:
    \begin{quote}\small
    At the time I originally discovered
the classification of diffeomorphism of surfaces, I was unfamiliar with
two bodies of mathematics which were quite relevant: first, Riemann surfaces,
quasiconformal maps and Teichm\"uller's theory; and second, Nielsen's theory
of the dynamical behavior of surface at infinity, and his near-understanding of
geodesic laminations.
\end{quote}

     In the same preface, Thurston writes: ``Dennis Sullivan first told me of some neglected articles by Nielsen which might be
relevant."     
     In a paper he wrote with M. Handel, titled \emph{New proofs of some results of Nielsen} \cite{HT1985}, Thurston  gave a new proof of his classification theorem using techniques from Nielsen's program. The relationship between the works of Thurston and those of Nielsen is also examined in the papers by Jane Gilman \cite{Gilman1981}, Richard Miller \cite{Miller} and Joan Birman \cite{Birman-Nielsen}.

           Talking about Nielsen, we are led to Nielsen's realization problem and the use of earthquakes in its solution.

           Earthquakes\index{earthquake} generalize the  Fenchel--Nielsen deformation operation  of cutting a hyperbolic surface along a simple closed geodesic and gluing it again after a twist. They are limits of sequences of such operations when the sequence of simple closed geodesics converges in Thurston's topology to a measured geodesic lamination. 
                       They became the canonical deformations between two hyperbolic structures after Thurston proved that for any two hyperbolic structures on a given surface, there is a unique left earthquake joining them. His proof is contained as an appendix in Steve Kerckhoff's paper \cite{K1}. They were the essential ingredient in Kerckhoff's proof of the Nielsen realization problem, which we review below.

Thurston wrote a paper on earthquakes on the hyperbolic plane \cite{Thurston-earthquakes}. In this paper, earthquakes are more naturally defined by cutting the hyperbolic plane along geodesics, taking limits of such operations, and studying the action on the unit circle at infinity. Considering this action on the universal covering  and on the circle at infinity solves the problems caused by the discontinuities of the map. At the same time, Thurston placed his theory in the setting of the universal Teichm\"uller space, the natural setting for deformations of the hyperbolic disc.  He obtained a new and elementary proof of the earthquake theorem.  He described this fact by the expression ``geology is transitive." In a set of notes he distributed in October 1987 on this new proof of the earthquake theorem, while he constructs the earthquake map using a homeomorphism of the circle at infinity of the hyperbolic plane, he writes:
\begin{quote}\small
This is closely connected to basic properties of convex hulls\index{convex hull} of sets in 3-space. Intuitively, imagine having disks with rubber arrows representing the identifications. Imagine some physical device which forces all the arrows to point counterclockwise: they bump against some barrier if you try to rotate them too far. You are allowed to move one of the disks by any isometry of the hyperbolic plane. You can kind of roll the disk around on the barriers through many different positions. This is very much like rocking a plane around on top of a wire which projects to a circle on a table. In the latter case, pushing straight down above one point finds the flat of the convex hull lying above a point inside the circle; in the former case, twisting at one point finds the stratum of the earthquake.

The earthquake theorem\index{earthquake theorem} can be proven by formalizing this argument.

\end{quote}

                             One of the first applications of Thurston's earthquake theorem was the proof of the Nielsen realization problem in 1980 by Kerckhoff \cite{K1, Kerckhoff}. The problem, formulated by Nielsen\index{Nielsen realization problem} in 1932,  asked whether any finite subgroup of the mapping class group of a surface can be realized as a group of homeomorphisms of this surface. In 1942,  Nielsen gave an affirmative answer in the case of finite cyclic groups.\footnote{Thurston and Handel note in their paper \cite{HT1985} that there should be a gap in Nielsen's proof of the fact that if a mapping class is periodic, then it contains a periodic homeomorphism of the surface. For the fact that Nielsen's proof is incorrect, they refer to Zieschang \cite{Zieschang}, and they declare that the known proofs of this fact use more sophisticated methods than those of Nielsen, e.g. actions on Teichm\"uller spce or Smith theory. Thus, they consider Nielsen's proof of the cyclic case as incomplete.} Fenchel extended Nielsen's result to the case of finite solvable groups.
                             There were several failed attempts to solve the Nielsen realization problem, namely, by Kravetz in 1959, based on the false assumption that the Teichm\"uller metric is negatively curved. But this failed proof had the advantage of putting the action of mapping classes on Teichm\"uller space at the center of the discussion.
                             Kerckhoff's proof is based on a convexity argument and a result of Thurston saying that any two points in Teichm\"uller space can be joined by a left earthquake.
                              In \S \ref{s:AdS} of this chapter, we shall talk about the  work of Geoffrey Mess in the early 1990s, who established a profound relation between earthquakes and the geometry of the convex core in  anti-de Sitter manifolds.    
                                                Besides the realization of finite subgroups, it was natural to address the same question for arbitrary subgroups.                                     Thurston asked the question of the lift of the whole group (see Problem 2.6 in Kirby's list \cite{Kirby}), and he conjectured that the answer is no. The conjecture was proved   by Markovic, for closed surfaces of genus $\geq 5$, in his paper \cite{M2007}, after Morita \cite{MoritaI} 
                                                obtained the same result for diffeomorphisms, using a more algebraic approach (cohomological obstructions). A year later, Markovic and \v Sari\' c completed the proof of Thurston's conjecture for the cases of genus 2 to 4 \cite{MS}.

         Talking about Thurston's work on Teichm\"uller space, we mention now his approach to the Weil--Petersson metric.\index{Weil--Petersson metric}
         
          Thurston introduced a Riemannian metric on Teichm\"uller space where the scalar product of two tangent vectors at a point represented by a hyperbolic surface is the second derivative of the length of a uniformly distributed sequence of closed geodesics on the surface. Thurston was motivated by the wish to have a metric defined using only the hyperbolic geometry of the surface, in contrast to the Teichm\"uller metric, which is based on quasiconformal theory, and to the Weil--Petersson metric, whose definition used the Petersson Hermitian product, defined in the context of modular forms and used by number theorists.
    
    Wolpert showed that Thurston's metric coincides with the Weil--Petersson metric \cite{Wolpert1986}. The consequence is that Thurston gave  a purely hyperbolic-geometric interpretation of the Weil--Petersson Riemannian metric on Teichm\"uller space.

       We now pass to Thurston's asymmetric metric.

In 1985--86, Thurston circulated a preprint titled  \emph{Minimal stretch maps between hyperbolic surfaces} \cite{Thurston1986} in which he 
introduced a non-symmetric metric on Teichm\"uller space which now bears the name \emph{Thurston metric}.\index{Thurston metric} The distance between two hyperbolic structures on a given surface is taken to be the logarithm of the infimum of the Lipschitz constants of homeomorphsims  that are homotopic to the identity,  where the distances used for the computation of the Lipschitz contant are the one of the first metric in the domain and the second metric in the range.  
 Thurston's paper is based on first principles (there is no appeal to any theory or any theorem except basic hyperbolic geometry).
 The paper was submitted to the journal \emph{Topology}. The referee sent a long report  asking for clarifications and references, and Thurston withdrew the paper.
In 1998, Thurston posted the article on the arXiv, and it remained unpublished. 

 Thurston's motivation was to develop a theory of Teichm\"uller space which is purely geometric, and which, like in his approach to the Weil--Petersson metric,  does not rely on quasiconformal mappings and quadratic differentials, but only on elementary hyperbolic geometry. He described a class of distinguished geodesics for this metric which he called stretch lines, he showed that any pair of points in the Teichm\"uller space can be joined by a concatenation of such lines, he showed that the metric is Finsler, and he described the dual unit sphere of the associated norm at each point of the cotangent space as an embedding  in this space of projective lamination space. 
  In the same paper, Thurston introduced his \emph{shear coordinates} for a surface decomposed into  ideal triangles, coordinates which have had an enormous impact, in the so-called higher Teichm\"uller theory and elsewehere, see e.g. \cite{LabM}; see also the generalization of these coordinates to the context of decorated Teichm\"uller theory \cite{Penner-Book}. 
Thurston's paper  \cite{Thurston1986} also contains the definition of an asymmetric norm on the tangent space to Teichm\"uller space which he called the \emph{earthquake norm}, which leads to another asymmetric metric on Teichm\"uller space. This norm is considered in the chapter by Barbot and Fillastre in the present volume \cite{Barbot-Fillastre}.

  In the first years after the preprint was released, little progress was made on this topic, one reason being that it took some time for the geometers to understand the ideas and the proofs contained in it. A survey of the results obtained in the first 20 years after Thurston's preprint was circulated, appeared in 2007, see \cite{PaT}. A set of open problems on Thurston's metric appeared in 2015 \cite{S}, after a conference held on this topic at the American Institute American Institute for mathematics in Palo Alto. We also mention the recent survey \cite{AP-Course}.

Thurston's metric led to the definition and study of analogous metrics in other settings:  surfaces with boundary \cite{2009f}, Euclidean structures on surfaces, \cite{BPT},  the Culler--Vogtmann outer space \cite{FM},
geometrically finite representations of fundamental groups of surfaces in   higher-dimensional  Lie groups  \cite{GK}, and there are many others developments, see e.g. \cite{AD, HP1}.

The next section, on fashion design, could have been included in the present one; it is also about Thurston's work on surfaces.

\subsection{Fashion design}
 
Let\index{fashion design} us start this section by  mentioning a  paper by Thurston and 
Kelly  Delp titled \emph{Playing with surfaces: Spheres, monkey pants, and zippergons} \cite{Delp}, one of the last papers of Thurston, written in 2011. In this paper, the authors describe a process, inspired by clothing design, of smoothing an octahedron to form a round sphere. They mention in the introduction several workshops and series of encounters they organized on devising schemes for designing pattern pieces to fit arbitrary shapes, including the human body. They declare: ``It was a very interesting but humbling experience, because
our initial assumption that familiar theoretical principles of differential geometry would do most of the work
was misleading."

Thurston was interested in the geometrical theory of clothes and the fitting of garments since his early work on surfaces. In the introduction to his paper \cite{Thurston1986}, he writes, after giving the definition of the best Lipschitz constant of maps in a given homotopy class:
\begin{quote}\small
This is closely related to the canonical problem that arises when a person
on the standard American diet digs into his or her wardrobe of a few years
earlier. The difference is that in the wardrobe problem, one does not really
care to know the value of the best Lipschitz constant---one is mainly
concerned that the Lipschitz constant not be significantly greater than 1. We
shall see that, just as cloth which is stretched tight develops stress wrinkles,
the least Lipschitz constant for a homeomorphism between two surfaces is
dictated by a certain geodesic lamination which is maximally stretched.
\end{quote}

 Thurston was not the first mathematician to think of the mathematical question of fitting a piece of fabric to some surface.  Pafnuty Chebyshev, back in the nineteenth century, thoroughly investigated the problem of fitting garments to a part of the human body. In particular, he thought about questions regarding the flexibility for a piece of fabric in order for it to approximate in the best exact form the part of body to which it is
designated and he established relations between this problem and other  mechanical problems he was studying, including the theory of linkages. He wrote an article on this topic, \cite{Chebyshev}.
A review of Chebyshev's work, making relations with his work and Euler's on geography and other problems related to surfaces,  is contained in the paper \cite{2016-Tchebyshev}.

  Thurston became eventually involved in fashion design. He worked with Dai Fujiwara, the creative director of the Japanese fashion designer Issey Miyake,  creating in 2010 a beautiful collection inspired by his eight geometries. On the occasion of the fashion show that took place at the Salon du Carrousel du Louvre in Paris, in March 2010, in which this collection was exhibited, 
   Thurston wrote a brief essay, distributed during the show, on beauty, mathematics and creativity.  
  Here is an excerpt:
 \begin{quote}\small
Many people think of mathematics as austere and self-contained. To the contrary, mathematics is a very rich and very human subject, an art that enables us to see and understand deep interconnections in the world. The best mathematics uses the whole mind, embraces human sensibility, and is not at all limited to the small portion of our brains that calculates and manipulates with symbols. Through pursuing beauty we find truth, and where we find truth, we discover incredible beauty.

The roots of creativity tap deep within to a place we all share, and I was
thrilled that Dai Fujiwara recognized the deep commonality underlying his
efforts and mine. Despite literally  and figuratively training and working on opposite ends of the earth, we had a wonderful exchange of ideas when he
visited me at Cornell. I feel both humbled and honored that he has taken up
the challenge to create beautiful clothing inspired by the beautiful theory
which is dear to my heart.
  \end{quote}

 In another article written on that occasion  for the fashion magazine \emph{Idom\'en\'ee}, Thurston made the following comment about the collection: 
\begin{quote}\small
The design team took these drawings as their starting theme and developed from there with their own vision and imagination. Of course it would have been foolish to attempt to literally illustrate the mathematical theory--- in this setting, it's neither possible nor desirable. What they attempted was to capture the underlying spirit and beauty. All I can say is that it resonated with me.
  \end{quote}

   Fashion design was for Thurston a ground for the combination of mathematics, art and practical applications, where the aesthetic component is pre-dominant.  In an interview released on the occasion of that fashion show, in which Thurston recounted how he came to contribute to the collection, he declared: ``Mathematics and design are both expressions of human creative spirit.''  About the aesthetical aspect of mathematics, Thurston had already written, in this 1990 \emph{Notices} article \cite{T-Math-edu}: 
  \begin{quote}\small
  My experience as a mathematician has convinced me that the aesthetic goals
and the utilitarian goals for mathematics turn out, in the end, to be quite close.
Our aesthetic instincts draw us to mathematics of a certain depth and connectivity.
The very depth and beauty of the patterns makes them likely to be manifested, in
unexpected ways, in other parts of mathematics, science, and the world.
  \end{quote}

                                \section{On Thurston's impact} 
                \subsection{The proof of the Smith conjecture}
  \label{Smith}
 
 The resolution of the Smith conjecture was the occasion for the  first major application of Thurston's uniformization theorem for Haken manifolds.\index{Smith conjecture}
 This conjecture says that if a cyclic group acts on $S^3$ by diffeomorphisms with one-dimensional fixed points, then it is topologically conjugate to the standard orthogonal action.
 The conjecture can be paraphrased as follows: a branched cyclic cover of a closed 3-manifold $M$ along a knot $K$ can be homeomorphic to the 3-sphere only if $M$ is also homeomorphic to the 3-sphere and $K$ is unknotted.
 
 The proof of the Smith conjecture was published as a book \cite{Sm}, which also contains a very comprehensive survey of Thurston's uniformization theorem written by Morgan.
The proof is divided into two cases depending on whether $M \setminus K$ contains a closed incompressible surface or not.
When it does, then results of Meeks--Yau and Gordon--Litherland,  contained in the book, show that the branched cover must also contain an incompressible surface, and hence cannot be the 3-sphere.
When it does not, then by the uniformization theorem,\index{uniformization theorem!Thurston} $M \setminus K$ is either a Seifert fibered manifold or hyperbolic.
In the former case, it is easy to see that the only possibility is that $M \setminus K$ is homeomorphic to $S^1 \times \mathbb R^2$.
In the latter case, an algebraic argument due to Bass implies that the holonomy representation of $\pi_1(M \setminus K)$ can be taken to lie in  a ring consisting of algebraic integers.
Then an argument of commutative algebra shows that the only possibility is that $K$ is unknotted.
 
After the publication of this book, Culler and Shalen \cite{CS} gave an alternative, more geometric approach to the algebraic argument in the last part of the proof.
They considered the algebraic set of characters of representations of the fundamental group of a hyperbolic manifold into $\mathrm{SL}_2 \mathbb C$, called the character variety.
\index{character variety}
Their alternative proof is obtained by considering the points at infinity of the character variety of $M \setminus K$, which gives a decomposition of $\pi_1(M \setminus K)$ by way of Bass--Serre theory.
This work of Culler--Shalen was generalized to a theory of compactification of character varieties\index{compactification!character variety} by Morgan--Shalen \cite{MS1, MS2, MS3}, which also gives an alternative proof of Thurston's compactness theorem for deformation spaces of acylindrical manifolds. 
                             
                             \subsection{The proofs of Ahlfors' conjecture, Marden's tameness conjecture, the ending lamination conjecture, and the density conjecture}
                             \label{Kleinian impact}
 The following four conjectures on Kleinian groups are contained in Thurston's list of  unsolved problems \cite{ThB}.
 \begin{enumerate}
 \item For any finitely generated Kleinian group, its limit set in the Riemann sphere either has measure $0$ or coincides with the entire sphere.
 This conjecture is originally due to Ahlfors \cite{Ahl}, and therefore called Ahlfors' conjecture.\index{Ahlfors' measure-0 conjecture}
 \item Any hyperbolic 3-manifold with finitely generated fundamental group is homeomorphic to the interior of a compact 3-manifold.
 This appeared first in Marden's paper \cite{Marden}.
 The property is called the topological tameness for the hyperbolic 3-manifold and also for the corresponding Kleinian group.\index{topological tameness}
The conjecture is called Marden's tameness conjecture.\index{Marden's tameness conjecture}
 \item
 If two hyperbolic 3-manifolds are homeomorphic and have the same parabolic locus, the same conformal structure at infinity and the same ending lamination then they are
 This is called the ending lamination conjecture.\index{ending lamination conjecture}
 \item
 Any finitely generated Kleinian group is an algebraic limit of geometrically finite Kleinian groups.
 This is called the (Bers--Sullivan--Thurston) density conjecture.\index{Bers--Sullivan--Thurston density conjecture} \index{density conjecture}
 \end{enumerate}     
 
The resolutions of these four conjectures proceeded in an intertwined way.
Thurston himself showed that algebraic limits of quasi-Fuchsian groups are geometrically tame.
Geometric tameness implies topological tameness, but is a stronger condition.
Thurston also showed that for geometrically tame Kleinian groups, Ahlfors' conjecture is true.
Bonahon \cite{Bo} clarified Thurston's notion of geometric tameness, and proved that any finitely generated Kleinian group that is not decomposed into a free product (i.e. any freely indecomposable Kleinian group) is geometrically and topologically tame, implying that Ahlfors' conjecture is also true for such Kleinian groups.

For freely decomposable Kleinian groups, Canary \cite{Ca} proved that topological tameness implies geometric tameness, and hence that Ahlfors' conjecture holds for topologically tame Kleinian groups.
Ohshika \cite{OhM} showed that any purely loxodromic algebraic limit of geometrically finite Kleinian groups is topologically tame unless the limit set is the entire sphere, and hence that Ahlfors' conjecture holds for any such algebraic limit.
Canary--Minsky \cite{CM} proved topological tameness for strong limits of topologically tame Kleinian groups.
Brock--Bromberg--Evans--Souto \cite{BBES} proved that every algebraic limit of geometrically finite Kleinian groups is topologically tame.
Finally, Agol \cite{Agp} and Calegari--Gabai \cite{CG} resolved Marden's tameness conjecture completely.

The ending lamination conjecture\index{ending lamination conjecture} was proved by Minsky \cite{MiT} for freely indecomposable Kleinian groups having a positive lower bound for the translation lengths (Kleinian groups\index{Kleinian group!bounded geometry} are then said to have bounded geometry).
Ohshika \cite{OhT} proved that the assumption of free indecomposability  can be removed still under the assumption of bounded geometry.
Minsky \cite{MiA} proved the ending lamination conjecture for once-punctured torus Kleinian groups.
The general ending lamination conjecture\index{ending lamination conjecture} was resolved by Minsky \cite{Mi} and Brock--Canary--Minsky \cite{BCM} using the work of Masur--Minsky \cite{MM1, MM2} on the geometry of curve complexes.
The proof relies on Thurston's idea of approximating the geometry of a neighborhood of an end by  pleated surfaces.
The point is that how pleated surfaces tend to the end is governed by a hierarchical structure of the curve complex, which was investigated in \cite{MM2}.

 The density conjecture was proved for Kleinian surface groups by Bromberg \cite{Br}, Brock--Bromberg \cite{BB} using Minsky's resolution of the ending lamination conjecture in the bounded geometry case.
 The general density conjecture was solved by Ohshika \cite{OhG} and Namazi--Souto \cite{NS} relying on the full resolution of the ending lamination conjecture.

 \subsection{The proof of the geometrization conjecture}
 The geometrization conjecture\index{geometrization conjecture} says that every closed  irreducible 3-manifold can be decomposed into geometric pieces by (Jaco--Shalen--Johannson) torus decomposition.
 Thurston's uniformization theorem says that this is true for Haken manifolds, but of course there are non-Haken manifolds: closed 3-manifolds with finite fundamental groups are clearly non-Haken, but there are also non-Haken manifolds with infinitely fundamental groups.
        
The geometrization\index{geometrization theorem} theorem was resolved by Perelman using the Ricci flow.\index{Ricci flow}
A Ricci flow is a deformation of Riemannian metric in the direction to reduce the variation of its Ricci curvature over the manifold, i.e. to average the Ricci curvature.
Hamilton  \cite{Ham} proved that any closed Riemannian 3-manifold with positive Ricci curvature is diffeomorphic to the $3$-sphere, making use of Ricci flows.
Perelman considered Ricci flows for general closed irreducible 3-manifolds. \cite{Perel1, Perel2, Perel3}
In contrast to Hamilton's case, the flow may encounter singularities.
Perelman showed that even in such cases, the deformation can be continued by rescaling and surgeries, and finally get to either constantly curved manifolds or Seifert fibrations.
This proves the geometrization conjecture, and in particular the Poincar\'{e} conjecture posed by Poincar\'{e} in 1904.

\subsection{The Waldhausen conjectures and the virtual fibering conjecture }\label{s:virtual}
In \cite{Waldpr}, Waldhausen\index{Waldhausen conjecture} posed the following two long-standing conjectures on 3-manifolds.
\begin{enumerate}
\item
The fundamental group of every closed irreducible 3-manifold either is finite or contains a closed surface group.
\item Every closed irreducible 3-manifold with infinite fundamental group is finitely covered by a Haken manifold.
This conjecture is now called the virtual-Haken conjecture.\index{virtual-Haken conjecture}
\end{enumerate}
The affirmative resolution of the second conjecture implies that of the first.

In his list of open questions in \cite{ThB}, Thurston took up the second conjecture again and added the following two stronger conjectures.

\begin{enumerate}
\setcounter{enumi}{2}
\item Every closed irreducible 3-manifold with infinite fundamental group is finitely covered by a 3-manifold with positive first Betti number.
\item Every closed hyperbolic 3-manifold is finitely covered by a surface bundle over the circle.
\end{enumerate}
Since it is easy to see that any Seifert fibered manifold with infinie fundamental group has a finite cover with positive Betti number, the affirmative resolution of (4) implies that of (3).

Conjecture (1) was known to hold for Seifert fibered manifolds. Therefore we have only to consider hyperbolic 3-manfolds.
The conjecture was proved in the case of arithmetic hyperbolic manifolds by Lackenby \cite{La}, and was
solved in general form by Kahn and Markovic \cite{KM2}.
Conjectures (2) and (4) (and hence also (3)) were solved by Agol, assuming the result of Kahn--Markovic, after partial results of Cooper--Long--Reid \cite{CLR} and Bergeron--Wise \cite{BW} among others.
These two works rely on quite different types of mathematics.

The resolution of (1) by Kahn and Markovic  took a rather analytic approach.
One first considers a pair of pants with geodesic boundaries in the hyperbolic manifold,  pastes a  pair of pants (with geodesic boundaries) to each of the boundary components, and then goes on pasting a pair of pants to each free boundary.
By a measure-theoretic argument, it is shown that after finitely many steps, the pair of pants comes back very close to the original one.
Pasting up all these pairs of pants, an immersed incompressible closed surface is obtained.

Agol's resolution of Conjectures (2) and (4) relies essentially on the study of CAT(0)-cube complexes\index{cube complex} started by Wise.
A cube complex is a complex made of finite-dimensional cubes, $[0,1]^n$, with isometric pasting maps.
A cube complex has a metric induced from the standard metrics on cubes, and is called\index{cube complex!CAT(0)} CAT(0) when it is non-positively curved in the sense of triangle comparison.
Bergeron--Wise \cite{BW} showed, using the work of Sageev \cite{Sa}, that the result of Kahn--Markovic cited above implies that the fundamental group of every closed hyperbolic 3-manifold acts freely and cocompactly on a CAT(0)-cube complex by isometries.
Haglund--Wise \cite{HW} considered \lq\lq hyperplanes'' in cube complexes, and introduced the notion of \lq\lq specialness'' for CAT(0)-cube complexes.
They then proved that if a hyperbolic group acts freely and cocompactly on a special CAT(0)-cube, then every quasi-convex subgroup is separable.
This implies that Conjecture (2) can be proved once we can show that every hyperbolic 3-manifold has a finite-sheeted cover whose fundamental group acts on a special CAT(0)-cube complex freely and cocompactly.
Agol proved that this is indeed the case.
Conjecture (4) was also resolved by combining this line of argument with previous work of Agol \cite{AgJ}.

  \subsection{The Ehrenpreis conjecture}
  The Ehrenpreis conjecture\index{Ehrenpreis conjecture} for Riemann surfaces  states that any two compact Riemann surfaces have finite sheeted unramified covers that are of the same genus and that are arbitrarily close to each other in the Teichm\"uller metric. It is not clear to the authors of this essay where and when exactly this conjecture was formulated for the first time. In their paper \cite{BN}, Biswas and Nag refer to it as an ``old conjecture which, we understand, is due to L. Ehrenpreis and C. L. Siegel."

              The conjecture was proved by J. Kahn
and V. Markovic in 2011 (the paper \cite{KM} was published in 2015). 
              We mention this here because the proof depends heavily on the geometric methods introduced by Thurston in the topology of surfaces and 3-manifolds. A crucial step in the proof is the construction of what the authors call a ``good" geodesic pair of pants decomposition of the surface, that is, a decomposition into pants whose cuff lengths are equal to some fixed large number.          
          Another major ingredient in the proof is an appeal to the proof of the surface subgroup  theorem and its proof  by the same authors,  which also makes heavy use of ideas inaugurated by  Thurston.     Thurston's influence on the subject is touched upon in the Kahn--Markovic ICM talk \cite{Markovic-ICM}.

              Sullivan and Thurston himself tried to prove this conjecture. In his approach to this question,  Sullivan introduced in the early 1990s \cite{Sulli} the notion of \emph{solenoid}, the inverse limit of the system of finite-sheeted branched covers of a fixed closed Riemann surface 
              of negative Euler characteristic. He introduced the Teichm\"uller space and the mapping class group of this object, and  studied their geometry and dynamics. The solenoid became an object of study in itself, see the reviews \cite{Penner-S, Saric}.

              In a memorial article on Thurston \cite{EMS-Surveys}, Sullivan recalls the following, from Milnor's 80th fest at Banff:  ``I recall a comment  whispered by Bill  who sat next to me during a talk by Jeremy Kahn about the Kahn--Markovich   proof  of the  subsurface conjecture\index{subsurface conjecture} from decades before.   Bill  whispered: `I missed the offset step'." (The ``offset step" referred to here is a step in the proof of  Kahn--Markovich which concerns the construction of pairs of pants with large cuffs).

\subsection{The Cannon--Thurston maps}
As an important step in the uniformization theorem for Haken manifolds, Thurston proved that any closed surface ($S$-)bundle over the circle with a pseudo-Anosov monodromy has a hyperbolic metric.
\index{Cannon--Thurston map}
His proof of this result uses the double limit theorem,\index{double limit theorem} which gives the Kleinian group  corresponding to the fiber.
The limit set of such a Kleinian surface group $G$ is the entire sphere, for it is a normal subgroup of a cocompact Kleinian group.
Cannon and Thurston \cite{CT} proved that there is a continuous map from $S^1$, which is the limit set of the Fuchsian group $\Gamma$ isomorphic to $\pi_1(S)$, to the limit set $S^2$ of $G$ which is equivariant under the action of $\Gamma$ on $S^1$ and $G$ on $S^2$.
Thus what they got is a $\pi_1(S)$-equivariant Peano map.

Thurston conjectured that such a map, called the Cannon--Thurston map, from the limit set of a convex cocompact Kleinian group to the limit set of isomorphic, possibly geometrically infinite group, which is invariant under the group action exists in general.
The existence of Cannon--Thurston maps was proved for Kleinian groups for freely indecomposable Kleinian groups with bounded geometry by Mahan in his paper \cite{Mahan1998}. The same result was obtained by Klarreich \cite{Kl}, based on the work of Minsky \cite{MiT}.
This was generalized to the case of punctured surface groups, still with the assumption of bounded geometry outside cusps, by Bowditch \cite{Bow}.
About ten years after he published his paper \cite{Mahan1998}, Mahan removed  an assumption he had there on the non-existence of parabolics \cite{Mahan2009}
McMullen, relying on \cite{MiA}, proved the existence of Cannon--Thurston maps for once-punctured torus Kleinian groups \cite{Mc}.
Finally, the general affirmative resolution of the conjecture was obtained by Manah Mj \cite{MjA, MjF} using the technique of \lq\lq electrocuting'' some parts of the manifold keeping the Gromov hyperbolicity.

Thurston also asked in the same list of unsolved problems if there is continuity  of movement of Cannon--Thurston maps with respect to the deformation of the Kleinian group.
Mahan--Series \cite{MjS1, MjS2} proved that when geometrically finite freely indecomposable Kleinian groups converge to a geometrically finite group algebraically, the Cannon--Thurston maps converge pointwise, but not necessarily uniformly.
They also proved that even if the limit is geometrically infinite, the uniform convergence is obtained provided that the limit is strong (i.e. it is both an algebraic and geometric limit).
In the case where the convergence is not strong, they gave an example when even the pointwise convergence fails.
Mahan--Ohshika \cite{MO} gave a necessary and sufficient condition for the pointwise convergence in the case where the sequence consists of quasi-Fuchsian groups.
\subsection{Anti-de Sitter geometry and transitional geometry}\label{s:AdS}

We start by recalling that for every $n\geq 2$, the $n$-dimensional Anti-de Sitter (AdS) space is a complete Lorentzian space of constant sectional curvature -1.   For $n\geq 2$, the model space of Anti-de Sitter space is the vector space $\mathbb{R}^{n+1}$ equipped with the  bilinear form of signaure $(n-1,2)$
\[<x,y>= -x_1y_1-x_2y_2+x_3y_3+\ldots+x_{n+1}y_{n+1}.\]

AdS space is the Lorentzian analogue\index{Lorentzian space} of hyperbolic space. 
In 1990, Geoffroy Mess  wrote a breakthrough paper\footnote{The paper was published in 2007 \cite{Mess2007};  see also the accompanying notes \cite{Andersson}.} in which he gave a completely new approach to Lorentzian geometry in dimension 2+1 and in which he proved a classification theorem for AdS spacetimes, that is, complete Lorentzian manifolds of constant negative curvature, obtained by taking a quotient of anti-de Sitter space by a discrete group of isometries acting freely.
 Mess's theory heavily uses techniques from Thurston's theory of low-dimensional geometry and topology, which he adapted to the Lorentzian setting. This includes a  Lorentzian version of the grafting\index{grafting!Lorentzian} operation for complex projective surfaces, the parametrization of the moduli space of complex projective structures as a bundle over the Teichm\"uller space of a surface whose fibers are measured laminations space, actions on $\mathbb{R}$-trees,  the holonomy map, convex hull and convex core constructions and a study of the geometry of the boundary of the convex core, bending and bending laminations,  representations of surface groups,  the analogue in AdS geometry of quasi-Fuchsian representations, the parametrization of moduli spaces of such representations  by two copies of Teichm\"uller space (an analogue of the Bers double uniformization theorem), and earthquakes,\index{earthquakes}  with a new proof of Thurston's earthquake theorem.

  The notions of Cauchy hypersurface\index{Cauchy hypersurface} and of globally hyperbolic\index{globally hyperbolic} AdS manifold turned out to be  central in this context:
A Cauchy hypersurface is a space-like hypersurface which intersects all inextendable time-like lines in the manifold in exactly
one point. An AdS manifold  (or more generally a  Lorentzian manifold) is said to be globally hyperbolic if it contains a Cauchy hypersurface.
 The notion of Cauchy surfaces was first introduced in the context of general relativity. 
 
 Mess, in his paper, gave a complete description of globally hyperbolic spacetimes of constant curvature
 with compact Cauchy surfaces in dimension 2+1. One of the results he obtained is the classification of proper isometric actions of discrete groups on Minkowski space.

  The introduction by Mess of Thurston's techniques in the setting of AdS geometry had an important impact on later research, and we mention some works done in this direction.

Bonsante and Schlenker, in their paper  \cite{BoS2009}  studied a space of adS manifolds with cone singularities,  and showed that this space is parametrized by the product of two copies of the
Teichm\"uller space of the surface with marked points (corresponding to the cone singularities). From this result they deduced an analogue of Thurston's theorem on the transitivity of  earthquakes for closed hyperbolic surfaces
with cone singularities with total angle less than  $\pi$. In the paper \cite{BoS2012}, the same authors   showed that it is possible to prescribe any two measured laminations filling a surface, to be the upper and lower measured bending laminations of the convex core of a globally hyperbolic AdS manifold, answering positively a question raised by Mess in his paper.
  
 Bonsante in his paper \cite{B2005} extended the study of globally hyperbolic flat spacetimes to higher dimensions. Among the tools he introduced is a notion of measured geodesic stratification which extends to higher dimensions that of measured geodesic lamination. Fillastre in \cite{Fillastre2011} studied Fuchsian polyhedra in such spaces, extending to this setting results of A. D. Alexandrov in \cite{Alexandrov}, Rivin--Hodgson in \cite{RH}, and  Labourie--Schlenker in \cite{LS2000} in which these authors  study convex Fuchsian surfaces in Lorentz spaces of constant curvature.

 In the paper \cite{BM2012}, Barbot and M\'erigot  established a relation between quasi-Fuchsian  and AdS representations which are Anosov in the sense of Labourie \cite{L2006}.

The work of Mess and later developments in AdS geometry are surveyed in Chapter 15 of the present volume, by Bonsante and Seppi \cite{BS-vol}, and in Chapter 16  by Barbot and Fillastre \cite{Barbot-Fillastre}.  The reader may also refer to the survey \cite{BeBo} by Benedetti and Bonsante and the works \cite{BMS1, BMS2, Barbot2018}.

 The list of open questions \cite{Questions-AdS} is another indication of the direction that this field took in the last few years, motivated by Thurston's ideas.

Talking about AdS geometry in relation with Thurston's work, we are led to transitional geometry, a topic also introduced by  Thurston.

A \emph{transition}\index{transitional geometry} between two geometries is a continuous path in the space of metrics on a manifold,   parametrized by an interval, say $(-1,1)$, where on the sub-interval $(-1,0)$ the manifold carries the first geometry (say hyperbolic), on the sub-interval $(0,1)$ it carries the other geometry (say AdS), and at $0$, the geometry is from a third type (say Euclidean).   
Thurston introduced the notion of transition\index{transitional geometry} between his eight 3-dimensional geometries in his proof of the orbifold theorem (see the  comments by Cooper, Hodgson and  Kerckhoff in \cite{CHK}),\footnote{The authors write in particular: ``Thurston outlined his proof on two occasions in courses at Princeton; in
1982 and again in 1984. On both occasions, due to running out of time,
the outline was incomplete in certain aspects at the end of the proof in the
collapsing case. In particular the Euclidean/spherical transition\index{transitional geometry} in the case
of vertices was treated in a few sentences."} and used the technique of the Ricci flow in this process. More details on this topic are given in \S \ref{geometrization} of the present article.

Following Thurston's ideas, there has been a recent activity in dimension $3$, on a continuous transition\index{transitional geometry} between the eight geometries, and also on varying continuously between Riemannian and Lorentzian geometries on orbifolds. 
Let us mention a few works on this subject.

Transitions between spherical and hyperbolic geometry, passing through Euclidean geometry, were studied by Cooper, Hodgson and  Kerckhoff in \cite{CHK}, Hodgson in \cite{Hodgson}, Boileau--Porti \cite{BP} and Porti in \cite{Porti2002}. 
In the last paper,  Porti investigated the appearance of orbifolds with geometry {\bf Nil} as
limits of rescaled hyperbolic cone manifolds. In his paper \cite{Porti2010}, he developed a 
theory of degeneration/regeneration between hyperbolic 2-orbifolds and hyperbolic cone 3-orbifolds.
In his paper with Weiss \cite{Porti-Weiss}, he developed a 
transition\index{transitional geometry} theory between Euclidean
cone manifolds and spherical or hyperbolic ones, with applications to questions of rigidity of Euclidean cone structures.
We also mention work on combinatorial transitions,\index{transitional geometry} by Kerckhoff and Storm in \cite{Kerckhoff-Storm}.

In the more recent paper \cite{CDW}, Cooper, Danciger and Wienhard  studied transitions between Thurston's geometries in the setting of projective geometry, 
They gave a complete classification of limits of
three-dimensional hyperbolic geometry inside hyperbolic geometry. They showed that the three Thurston geometries 
$\mathbb{E}^3$, {\bf Nil,} and {\bf Sol} appear among the limits, but the other Thurston geometries do not.

 In the papers \cite{Danciger2013, Danciger2014}, Danciger studied a smooth transition\index{transitional geometry} between the hyperbolic and AdS  geometry of 3-manifolds, passing through a transversely hyperbolic 1-dimensional foliation of the manifold.  
In particular, in the first paper, he introduced, in a study of the transition geometry between hyperbolic 
and AdS geometry, a transitional projective geometry he called half-pipe geometry.

Two-dimensional transitional geometry,\index{transitional geometry} from a completely different point of view, based on the notion of ``coherent geometry", is studied by  A'Campo and Papadopoulos \cite[Chapter 9]{2011f} and \cite{ACP1}.

                               \subsection{Linkages}  The theory of linkages\index{linkage} is a classical subject that combines topology, real-analytic geometry and  mechanical constructions. It is not surprising that this topic attracted Thurston's attention. Although he wrote very little on it, he influenced the works of several authors, in private conversations and in lectures. Among them,  we mention  Henry King, Misha Kapovich, 
John Millson and Alexei Sossinsky.  
                                                     
                  In the survey  \cite{King} on planar linkages,  King recalls  that he first heard of the subject in a talk by Thurston at the Institute for Advanced Study, in the
mid-1970s. He remembers that Thurston gave a proof of what is called now the\index{Thurston signature theorem} \emph{Thurston signature theorem.}\footnote{The word ``signature" refers here to a person signing her name; this is not to be confused with notions like Rokhlin's signature of a manifold.} The statement can take different forms, one of them being that for any signature, one can construct a planar linkage  
                   that approximates it arbitrarily closely. In other words, one can find a linkage such that the locus of one of its vertices (or of a set of vertices), when it runs through all its possible positions, is arbitrarily close to the signature.  
                   Another (related) result of  Thurston mentioned in King's paper concerns the realization of any compact smooth manifold as a configuration space of
a linkage.
King writes:  
\begin{quote}\small
 As far as I can tell, Thurston never wrote these results up, so [the Thurston signature theorem] must remain vague.
Occasionally since then I have been contacted by an engineer interested in these results,
but I could not recall anything about Thurston's proof so I could not help them. Then
       recently, Millson started asking me lots of questions on real algebraic sets. He and
 Kapovich were writing up proofs of the results [of Thurston] above. In the course of doing
so, they discovered and solved some problems overlooked by previous literature.
\end{quote}
            
             Kapovich and Millson, in the paper they wrote on the subject \cite{KaM},  make a reference to the work of the 19th century mathematician Alfred Kempe  \cite{Kempe} who studied linkages  and obtained weak forms of  Thurston's results. This  fact is recurrent in Thurston's work: He used to develop theories from scratch, and it happened that he realized that some of his ideas or results were discovered by others, in general, several years, and sometimes decades before him. In fact, Thurston was reviving classical subjects.  We already mentioned such instances in the section concerning his work on surface diffeomorphisms (\S \ref{s:surfaces}) and in that on circle packings (\S \ref{s:circle}). There are many other examples.               
                  
               Kapovich and Millson write in their paper: 
               \begin{quote}\small
               The first precise formulation of a theorem of the above type was given by W. Thurston who stated a
version of Corollary C about 20 years ago and has given lectures on it since. He realized that
such a theorem would follow by combining the 19th century work on linkages (i.e. Kempe's
theorem) with the work of Seifert, Nash, Palais and Tognoli. However, Thurston did not write
up a proof so we have no way of knowing whether he overcame the problems discussed above in
the 19th century work on linkages. There is also ambiguity concerning which theorem Thurston
formulated in his lectures, we heard three different versions from three sources.
\end{quote}
                                                                     
                       Sossinsky writes in his survey  \cite{Sossinsky} that he was introduced to the theory of linkages by Alexander Kirillov, after the latter returned to Moscow  from a stay in the US during which he had heard one of Thurston's talks on the subject. In turn, Sossinsky introduced the subject to several Russian mathematicians who started working on it. He writes that Thurston was mainly interested in the topology of configuration spaces of planar linkages, and that he considered two types of problems, which were the main problems in the field: (i) the so-called 
direct problem (``configuration"): given a planar linkage (or a class of planar
linkages), to ask for a description of the corresponding configuration space(s)  ; (ii) the
inverse problem (``universality"): given a topological space or an algebraic variety (or a class of such objects), to find a planar linkage whose configuration space
is this space (or a class of linkages whose configuration spaces are in the given class). On problem (ii), Sossinsky mentions a version of Thurston's 
 signature theorem saying that for any real-algebraic  curve in the plane there exists 
 a planar linkage which draws it. Sossinsky's article \cite{Sossinsky} contains a beautiful historical introduction to the subject.  
 
  In a post on MathOverflow, Kevin Walker recalls that as an undergraduate student of Thurston at Princeton, the latter told  him about the strategy of the proof of the signature theorem.  This proof included a use of Nash's theorem saying that any smooth manifold is diffeomorphic to a real algebraic set, which reduces the problem to that of devising planar linkages implementing addition and multiplication of real numbers and showing how to combine these linkages.  Walker wrote his bachelor thesis on linkages. He posed there a conjecture, about recovering the relative lengths of the bars of a linkage
from intrinsic algebraic properties of the cohomology algebra of its configuration space. The conjecture was proved several years later by Farber, Hausmann and Sch\"utz \cite{FHS}.
              
 Thurston's  popular science article with Jeff Weeks \cite{TWeeks}, published in 1984, contains several passages on linkages. In particular, we find there the description of a simply defined linkage (which 
 was studied later by several authors under the name Thurston--Weeks triple linkage)\index{Thurston--Weeks triple linkage} whose configuration space has a very interesting
topology.

                        In a correspondence with the second author of this chapter, Bill Abikoff wrote: ``Thurston was characteristically terse in his discussion of spaces formed by flexible linkages.  His response to the question of which topological spaces appear as the configuration space of a flexible linkage was: \emph{all}."                               
 \subsection{Higher Teichm\"uller theory}

Thurston is a forerunner of higher Teichm\"uller theory.\index{higher Teichm\"uller theory}\index{Teichm\"uller theory!higher} He was the first to emphasize the importance of the study of connected components of the representation variety of the fundamental groups surfaces into Lie groups other than the group $\mathrm{PSL}(2,\mathbb{R})$.  He was also the first to revive ideas of Ehresmann from the mid 1930s, highlighting the holonomy as a map
from the deformation space of geometric structures structures into the representation variety, making this a general guiding principle for the classification
of locally homogeneous structures.  We refer the reader to Goldman's article \cite{Goldman} in which he talks about  an \emph{Ehresmann--Weil--Thurston holonomy principle}.\index{holonomy principle} 
Labourie and McShane use the expression ``Higher  Teichm\"uller--Thurston theory" for the study of a specific component of the representation space of a surface group of genus  in $\mathrm{PSL}(n,\mathbb{R})$. In their paper \cite{LabM},   they extend Thurston's shear coordinates to the  setting of Hitchin representations  of fundamental groups of surfaces and they prove a McShane--Mirzakhani\index{McShane---Mirzakhani identity}
identity in that setting. Vlamis and Yarmola use the same expression in the paper where they prove a Basmajian identity in higher  Teichm\"uller--Thurston theory \cite{VY}.
Among the large number of  results in higher Teichm\"uller theory that are inspired by Thurston's work on surfaces, we mention Labourie's work on representations of surface fundamental groups into $\mathrm{PSL}(n,\mathbb{R})$, and in particular his discovery of a curve which is the limit set of the quasi-Fuchsian representation in this setting \cite{L2006}. We also mention  
the generalization of Thurston's shear coordinates to the context of  decorated representations
into split real Lie groups by Fock and Goncharov \cite{FG}, and
 the generalization of 
Thurston's compactification of Teichm\"uller space to compactifications of spaces of various sorts of representations of finitely generated groups  (see e.g. \cite{Parreau} for representations into reductive Lie groups). One should also mention the recasting of Thurston's compactification of Teichm\"uller space from the  point of view of the character variety, in terms of  group actions on $\Lambda$-trees, by Morgan and Shalen, see \cite{MS1}. Finally, we mention the work done on the pressure metric on higher Teichm\"uller spaces (in particular for Anosov representation)
a higher-generalization of Thurston's version of the Weil--Petersson metric on Teichm\"uller space, see \cite{BCLS, BCS}.

  \subsection{The Grothendieck--Thurston theory}

Alexander Grothendieck,\index{Grothendieck--Thurston theory} at several places of his manuscript  \emph{Esquisse d'un programme} \cite{Esquisse} (released in 1984), in which he  introduced the theory of \emph{dessins d'enfants}\index{dessin d'enfants}  and where he set out the basis of the theory that later on  became known as Grothendieck--Teichm\"uller theory, mentions Thurston's work as a source of inspiration.   On p. 12 of his manuscript, Grothendieck writes: ``The lego-Teichm\"uller toy which I am trying to describe, arising from motivations and reflections of absolute algebraic geometry on the field $\mathbb{Q}$, is very close to Thurston's hyperbolic geodesic surgery."  Grothendieck drew a parallel between his own algebraic constructions in the field $\mathbb{Q}$ of rational numbers and what he calls Thurston's ``hyperbolic geodesic surgery" of a surface by pairs of pants decompositions. He outlined in this paper a principle which today bears the name ``Grothendieck\index{Grothendieck reconstruction principle} reconstruction principle," or  the ``two-level principle."\index{two-level principle} In broad terms, the principle says that some important geometric,
algebraic and topological objects that are associated with a surface $S$ (e.g.
the Teichm\"uller space, the mapping class group, the space of measured foliations, spaces of representations of its fundamental group, etc.) can be reconstructed
from the ``small" corresponding spaces associated with the (generally infinite) set of level-zero,
level-one and level-two essential subsurfaces of $S$. Here, the ``level" of a surface is the number of simple closed curves that are needed to decompose it into pairs of pants. Thus, level-zero surfaces are pairs of pants, level-one surfaces are tori with one hole or spheres with 4 holes, level-two surfaces are the 2-holed tori and 5-holed spheres etc.
 The geometric structures on the
level-zero spaces are the building blocks of the general structures, and the structures
on the level-one and the level-two spaces are the objects that encode the gluing. There is a group-theoretic flavor where the level-one surfaces play the role of generators and the level-two surfaces are the corresponding relators.
 Paraphrasing
Grothendieck from his Esquisse d'un programme, ``the Teichm\"uller tower\index{Teichm\"uller tower}
can be reconstructed from level zero to level two, and in this reconstruction, level-one
gives a complete set of generators and level-two gives a complete set of relations."
Grothendieck made a comparison with analogous situations in algebraic geometry,
in particular in reductive group theory, where the semi-simple
rank of a reductive group plays the role of ``level."

The reconstruction principle was used (without the name) in the paper by Hatcher and Thurston, \emph{A presentation for the mapping class group of a closed orientable surface} \cite{HT}, published in 1980. In this paper, the authors find a presentation of the mapping class group in which the generators and the relations,
which correspond to moves in the
pants decomposition complex, are all supported on the level-two surfaces of the given
topological surface. The reconstruction principle\index{reconstruction principle} appears in the same paper at the level of functions: the authors use Cerf theory (the study singularities in the space of
smooth functions on the surface) in a construction which is also limited to the level-one and level-two subsurfaces.  
The analogy between Grothendieck's and Thurston's theories is expanded  in Feng Luo's paper \emph{Grothendieck's reconstruction principle and 2-dimensional topology and geometry} \cite{Luo}.   
 
   On p. 41 of his manuscript \cite{Esquisse},  Grothendieck  formulates and discusses a conjecture concerning the canonical realization of conformal structures on surfaces by complex algebraic curves. He then declares: ``An elementary familiarization with Thurston's beautiful ideas on the construction of Teichm\"uller space in terms of a very simple game of Riemannian hyperbolic surgery reinforces my premonition."

 Grothendieck also used ideas of Thurston in his works on the actions of the absolute Galois group and in profinite constructions in Teichm\"uller's theory. The author may refer to the surveys \cite{AJP, Uludag}.
 The same ideas are also developed in his manuscript \emph{Longue marche \`a travers la th\'eorie de Galois} \cite{Longue}, written around the same period.  
 
 At the University of Montpellier, where he worked for the last 15 years of his academical life, Grothendieck conducted a seminar on Thurston's theory on surfaces. 
 
 Grothendieck again mentions Thurston's work on surfaces in his mathematical autobiography, \emph{R\'ecoltes et semailles} \cite[\S 6.1]{RS}.  In that manuscript he singles out twelve themes that dominate his work and which he describes as ``great ideas" (grandes id\'ees).  Among the two themes he considers as being the  most important is what he calls the ``Galois--Teichm\"uller yoga",\index{Galois--Teichm\"uller yoga} which is precisely the topic that now bears the name  Grothendieck--Teichm\"uller theory\index{Grothendieck--Teichm\"uller theory} \cite[\S 2.8, Note 23]{RS}.
    
     Grothendieck and Thurston had different approaches to Teichm\"uller space, because the motivations were different (algebraic geometry and low-dimensional topology), but reuniting the two approaches is still now a challenging  field of research. Mapping class groups of surfaces occur in the Grothendieck setting in the form of the so-called  Grothendieck--Teichm\"uller group\index{Grothendieck--Teichm\"uller group} and in the Teichm\"uller tower,\index{Teichm\"uller tower} built out of finite type surfaces. The curve complexes and other simplicial complexes from low-dimensional topology have their analogues in this theory. In fact, some of the tools in Grothendick's theory are profinite versions of notions discovered by Thurston. The interested reader may refer to the surveys \cite{Funar, Uludag}. Conversely, Grothendieck's dessins d'enfants were studied by several authors in the setting of Thurston's theory; we refer to the surveys \cite{Herrlich, Harvey}.
          
 To close this section, let us recall that both Grothendieck and Thurston campaigned against military funding of mathematics. In France, Grothendieck resigned abruptly  from his position at IH\'ES after he learned that the institute was run partially by military funds. Ten years later, in the US, Thurston was thoroughly involved in a campaign against military funding of mathematics. He wrote several letters to the editors of the \emph{Notices of the AMS}, see e.g. his article   \emph{Military funding in mathematics} \cite{T-Funding}.

\section{In guise of a conclusion}
A description of Thurston's work would be incomplete without a few words on his personality. 
 
 Thurston  valued the notion of mathematical community, and he was pleased to see that he could share his ideas with more and more people. Beyond mathematics, his militancy for a good educational system, for the protection of nature and for a clean environment, his search for beauty, his gentleness, his humbleness, his honesty, and his care for people around him and for humanity in general  were exceptionally high. He was a rebel in every sense of the word.

%  \vfill\eject        
%\-\vfill
%\thispagestyle{empty}
%
%
% \includegraphics[width=1\linewidth]{Thurston2.jpg}   
%  
% \bigskip
%        
%\-\vfill
%\thispagestyle{empty}
% 
%
% \includegraphics[width=1\linewidth]{Thurston1.jpg}    
% 
% 
%  \caption{\small  Clay conference in Paris, Oceanographic Institute, June 2010. @ Atelier EcoutezVoir. 
%
%
%\vfill\eject
%

    \end{document}